\documentclass[11pt,twoside,leqno,a4paper]{article}

\usepackage{geometry, amsfonts, amsthm, amsmath, indentfirst, graphicx, enumerate,xcolor, subcaption, lineno, cite, float}
\usepackage{paralist}
\usepackage{amssymb}
\usepackage{setspace} % I use it to reduce the space between the lines in the table of contents

\usepackage{extarrows}% for making wide '=' symbols. See the command \xlongequal

\usepackage[nottoc,notlot,notlof]{tocbibind}%for appearing references in the list of contents

\RequirePackage[breaklinks,colorlinks]{hyperref}
\hypersetup{linkcolor=blue,citecolor=blue,urlcolor=blue}

\usepackage{titlesec}
\titleformat{\section}{\filcenter\Large\bfseries}{\normalfont\thesection.}{0.5em}{}
\titleformat{\subsection}{\filcenter\bfseries}{\normalfont\thesubsection.}{.5em}{}

\newtheorem{mythm}{\,\,\,\,\,\, \normalfont\scshape Theorem}[section]
\newtheorem{theorem}{\,\,\,\,\,\, \normalfont\scshape Theorem}
\newtheorem{mycor}[mythm]{\,\,\,\,\,\, \normalfont\scshape Corollary}
\newtheorem{mylem}[mythm]{\,\,\,\,\,\, \normalfont\scshape Lemma}
\newtheorem{myprop}[mythm]{\,\,\,\,\,\, \normalfont\scshape Proposition}

\newtheorem{myrem}[mythm]{\,\,\,\,\,\, \normalfont\it Remark}
\newtheorem{mynotation}[mythm]{\,\,\,\,\,\, \normalfont\it Notation}

\newtheorem{myassumption}{\,\,\,\,\,\, \normalfont\it Assumption}

\newtheorem{theoremA}{\,\,\,\,\,\, \normalfont\scshape Theorem}

\newtheorem{theoremB}{\,\,\,\,\,\, \normalfont\scshape Theorem}

\definecolor{mygreen}{rgb}{0.50,0,0.50}

\numberwithin{equation}{section}

\usepackage{etoolbox}
\makeatletter
\patchcmd{\endmyassumption}{\@endpefalse}{}{}{}
\makeatother

\usepackage{authblk}

\AtEndDocument{%
\par
\medskip
\begin{tabular}{@{}l@{}}%
{\small Instituto de Ci\^encias Matem\'aticas e Computac\~ao (ICMC), Universidade de S\~ao Paulo, S\~ao Carlos, Brazil}\\
\textit{\small E-mail address}: \href{mailto:sajjad.bakrani@icmc.usp.br}{\texttt{\small sajjad.bakrani@icmc.usp.br}}
\end{tabular}}

\begin{document}

\title{\bf Dynamics near homoclinic orbits to a saddle in four-dimensional systems with a first integral and a discrete symmetry}
\author{\textsc{Sajjad Bakrani}}
\date{November 3, 2024}
\maketitle

\begin{abstract}
We consider a $\mathbb{Z}_{2}$-equivariant 4-dimensional system of ODEs with a smooth first integral $H$ and a saddle equilibrium state $O$. We assume that there exists a transverse homoclinic orbit $\Gamma$ to $O$ that approaches $O$ along the nonleading directions. Suppose $H(O) = c$. In \cite{Bakrani2022JDE}, the dynamics near $\Gamma$ in the level set $H^{-1}(c)$ was described. In particular, some criteria for the existence of the stable and unstable invariant manifolds of $\Gamma$ were given. In the current paper, we describe the dynamics near $\Gamma$ in the level set $H^{-1}(h)$ for $h\neq c$ close to $c$. We prove that when $h < c$, there exists a unique saddle periodic orbit in each level set $H^{-1}(h)$, and the forward (resp. backward) orbit of any point off the stable (resp. unstable) invariant manifold of this periodic orbit leaves a small neighborhood of $\Gamma$. We further show that when $h > c$, the forward and backward orbits of any point in $H^{-1}(h)$ near $\Gamma$ leave a small neighborhood of $\Gamma$. We also prove analogous results for the scenario where two transverse homoclinics to $O$ (homoclinic figure-eight) exist. The results of this paper, together with \cite{Bakrani2022JDE}, give a full description of the dynamics in a small open neighborhood of $\Gamma$ (and a small open neighborhood of a homoclinic figure-eight).
\end{abstract}

%\blfootnote{\keywords{homoclinic, super-homoclinic, invariant manifold, coupled Schr\"odinger equations}}

%{\small\tableofcontents}

\begin{spacing}{0}
\small\tableofcontents
\end{spacing}

\addtocontents{toc}{\protect\setcounter{tocdepth}{1}}
% This command controls what appears on the table of contents.

%\listoffigures

\section{Introduction}

This paper aims to describe the dynamics near homoclinic orbits in 4-dimensional Hamiltonian systems (or, more generally, systems with smooth first integrals). An orbit $\Gamma$ of a system of ODEs is called homoclinic if it lies at the intersection of the stable and unstable manifolds of an equilibrium point $O$ of the system. As a result, a homoclinic trajectory $\Gamma = \{x\left(t\right)\}$, where $x$ and $t$ are the state and time variables, respectively, converges to $O$ as $t\rightarrow -\infty$ and $t\rightarrow +\infty$. Having a Hamiltonian (or more generally, a smooth first integral) $H$, we have that both the homoclinic orbit $\Gamma$ and the equilibrium $O$ are in the same level set of $H$. To describe the dynamics in a small neighborhood $U$ of $\Gamma \cup \{O\}$, it is sufficient to study the dynamics in the restriction of $U$ to the level sets $\{H = h\}$ for $h$ close to $H(O)$. This is because each level set of $H$ is invariant with respect to the flow of the system\footnote{Dynamics near homoclinic orbits in the non-conservative case has been studied in the literature too, see e.g. \cite{Dimathesis, Homburg1996memoirs, HomburgSandstede2010, Sandstede2000center, Shashkov1999existence, Shilnikov1965}.}.

Generically, the homoclinic orbit $\Gamma$ approaches $O$ for both $t\rightarrow -\infty$ and $t\rightarrow +\infty$ along the leading directions (the directions associated with those eigenvalues of the linearized system at $O$ that are nearest to the imaginary axis). Dynamics near homoclinic orbits in this generic scenario has been studied in the literature. In the case that the leading eigenvalues of the linear part of the system at $O$ are real, Turaev and Shilnikov \cite{DimaShilnikov1989, Turaev2014} (see also \cite{ShilnikovDimaSuperhomoclinic1997}) proved that the only orbits that remain in the level set $\{H = H(O)\}$ are the homoclinics and the equilibrium $O$, while the set of the orbits that remain in the level set $\{H = h\}$ for $h\neq H(O)$ sufficiently close to $H(O)$ is a hyperbolic set (the restriction of the dynamics to this set is topologically equivalent to the suspension over a Markov chain). The dynamics near homoclinic orbits when the leading eigenvalues are non-real has also been studied in the literature, and it is shown that it can be very nontrivial \cite{Devaney1976, BelyakovShilnikov1990, Lerman1991, Lerman2000, Lerman1997, Buffoni1996, Barrientos2016robust}.

This paper considers a degenerate scenario where the homoclinic orbits approach a saddle equilibrium with real leading eigenvalues along nonleading directions as $t\rightarrow \pm\infty$. This scenario can become generic in the presence of symmetries as the existence of symmetries implies the existence of (low-dimensional) invariant vector subspaces, and one may expect that stable and unstable manifolds of the equilibrium are more likely to intersect when the system is reduced to the invariant subspace. In this paper, we also consider the simplest case of a $\mathbb{Z}_{2}$-symmetry. The dynamics in the level set $\{H = H(O)\}$ for this setting was studied in \cite{Bakrani2022JDE}. It was established there that the presence of the $\mathbb{Z}_{2}$-symmetry may imply the existence of stable and unstable invariant manifolds of the homoclinic orbits, which may lead to the emergence of superhomoclinic orbits (orbits that lie in the intersection of the stable and unstable manifolds of the homoclinic orbits), and consequently, infinitely many multi-pulse homoclinic loops.

Here, we study the dynamics near the homoclinic orbits in the level sets $\{H = h\}$ for $h \neq H(O)$ sufficiently close to $H(O)$. We prove that the set of the orbits in these level sets that stay close to the homoclinic orbits for forward (resp. backward) time is either empty or constitutes of some periodic orbits (which are saddle in their level sets) and their stable (resp. unstable) invariant manifolds. This result, together with the one in \cite{Bakrani2022JDE}, describes the dynamics in a small open neighborhood $U\subset \mathbb{R}^{4}$ of the homoclinic orbits.

One example of our setting is the Coupled nonlinear Schr\"odinger equations (CNLSE) given by
\begin{equation*}
\begin{aligned}
i\Psi_{t} + \Psi_{xx} + 2\left(\alpha \lvert\Psi\rvert^{2} + \lvert\Phi\rvert^{2}\right)\Psi &= 0,\\
i\Phi_{t} + \Phi_{xx} + 2\left(\lvert\Psi\rvert^{2} + \beta \lvert\Phi\rvert^{2}\right)\Phi &= 0,
\end{aligned}
\end{equation*}
where $\Psi$ and $\Phi$ are complex-valued functions of $(t,x)$, and $\alpha$ and $\beta$ are complex constants. The CNLSE is considered as a basic model for light propagation. It also appears as a universal model of behavior near thresholds of instabilities (see e.g. \cite{Kirrmann_CNLSE_1992}). Taking $\alpha$ and $\beta$ to be positive real constants and looking at steady state solutions of the form
\begin{equation*}
\Psi(t,x) = e^{i\omega_{1}^{2} t} \psi(x) \qquad \mathrm{and}\qquad \Phi(t,x) = e^{i\omega_{2}^{2} t} \phi(x),
\end{equation*}
the CNLSE can be written in the form of a 4-dimensional Hamiltonian system with two $\mathbb{Z}_2$-symmetries and a pair of homoclinic figure-eights (see e.g. \cite{BakraniPhDthesis}). Our results, together with the ones in \cite{Bakrani2022JDE} and \cite{Turaev2014}, describe the dynamics near these homoclinic orbits.

\subsection{Problem setting and results}\label{Setting and results}

\begin{figure}
\centering
\includegraphics[scale=.18]{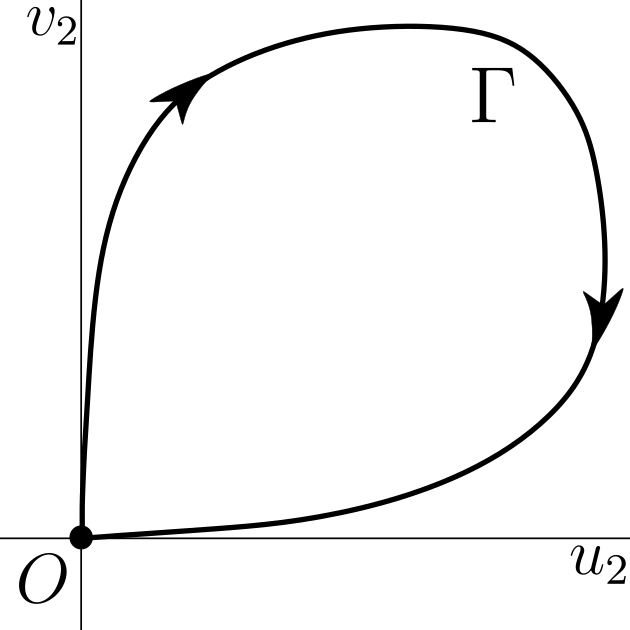}
\caption{\small The transverse homoclinic loop \(\Gamma\) in the invariant plane \(\lbrace u_{1}=v_{1}=0 \rbrace\).}
\label{Figure8o98b87rv8i6}
\end{figure}

Consider a \(\mathcal{C}^{\infty}\)-smooth 4-dimensional system of differential equations
\begin{equation}\label{eq100}
\dot{x}=X(x), \quad x\in \mathbb{R}^{4},
\end{equation}
with a hyperbolic equilibrium state $O$ at the origin. We assume that this system possesses a \(\mathcal{C}^{\infty}\)-smooth first integral\footnote{In some literature, it is called 'energy function', 'integral of motion' or 'constant of motion'.} \(H: \mathbb{R}^{4}\rightarrow \mathbb{R}\), i.e. $H$ is constant along any orbit of system (\ref{eq100}) and its restriction to any open subset of $\mathbb{R}^4$ is a non-constant function. One may consider Hamiltonian systems as a natural example of such a setting, however, the symplectic structure is not required here. Since $H$ is constant along orbits of (\ref{eq100}), we have $H^{\prime}(x) \cdot X(x) \equiv 0$. In particular, \(H^{\prime}(O) X^{\prime}(O) \equiv 0\), and since $O$ is hyperbolic, this implies that the linear part of \(H\) at \(O\) vanishes. We consider a generic scenario in which the quadratic part of \(H\) at \(O\) is nondegenerate\footnote{A quadratic form in $n$ variables $x_1, \ldots, x_n$ is nondegenerate if it can be written as $(x_1, \ldots, x_n) A (x_1, \ldots, x_n)^{\top}$ for a nonsingular symmetric $A\in\mathbb{R}^{n\times n}$.}. Then, there exists a linear change of coordinates which reduces system (\ref{eq100}) near \(O\) to the form
\begin{equation}\label{eq300}
\dot{u}=-Au + o\left(\lvert u\rvert, \lvert v\rvert\right), \qquad \dot{v}=A^{T}v + o\left(\lvert u\rvert, \lvert v\rvert\right),
\end{equation}
where \(u=\left(u_{1},u_{2}\right)\in \mathbb{R}^{2}\), \(v=\left(v_{1},v_{2}\right) \in \mathbb{R}^{2}\) and \(A\in\mathbb{R}^{2\times 2}\) is a matrix whose eigenvalues have positive real parts. Moreover, the first integral $H$ is brought to the form
\begin{equation}\label{eq400}
H=\langle v, Au\rangle + o\left(\|u\|^{2} + \|v\|^2\right),
\end{equation}
where \(\langle\cdot, \cdot\rangle\) is the standard inner product on \(\mathbb{R}^{2}\) \cite{BakraniPhDthesis}.
\begin{myassumption}\label{assumption30}
System (\ref{eq300}) is invariant with respect to the symmetry
\begin{equation}\label{eq425}
(u_{1}, v_{1})\leftrightarrow(-u_{1}, -v_{1}).
\end{equation}
\end{myassumption}
Following this assumption, the plane \(\lbrace u_{1}=v_{1}=0 \rbrace\) is invariant with respect to the flow of system (\ref{eq300}). On the other hand, the 2-dimensional stable and unstable invariant manifolds of the saddle $O$, \(W^{s}\left(O\right)\) and \(W^{u}\left(O\right)\), respectively, lie in the 3-dimensional level set \(\lbrace H = 0 \rbrace\), and so one may expect that they intersect transversely in that level, producing a number of homoclinic loops. We consider the following specific case.
\begin{myassumption}\label{assumption50}
There exists a homoclinic orbit \(\Gamma\) of the transverse intersection of \(W^{s}\left(O\right)\) and \(W^{u}\left(O\right)\) lying in the invariant plane \(\lbrace u_{1}=v_{1}=0 \rbrace\) (see Figure \ref{Figure8o98b87rv8i6}).
\end{myassumption}

The action of the symmetry commutes with the linear part of system (\ref{eq300}) implying that \(A\) is a diagonal matrix, i.e. $A = \mathrm{diag}(\lambda_1, \lambda_2)$ for some real $\lambda_1, \lambda_2 > 0$. Without loss of generality, assume $\lambda_{2} \geq \lambda_{1}$. Taking this into account and by rescaling (\ref{eq400}), we write the first integral $H$ as
\begin{equation}\label{First_integral}
H = \gamma u_{1}v_{1} - u_{2} v_{2} + o\left(\|(u, v)\|^2\right),\qquad \mathrm{where}\,\,\, \gamma := \lambda_{1} \lambda_{2}^{-1} \leq 1.
\end{equation}

This paper studies the dynamics in a sufficiently small neighborhood of the homoclinic loop $\Gamma \cup \{O\}$ (and the homoclinic figure-eight discussed later). To this end, we first reduce system (\ref{eq300}) near the equilibrium $O$ to a normal form; it is done for the cases $\lambda_{1} = \lambda_{2}$, $\lambda_{1} < \lambda_{2} \leq 2\lambda_{1}$ and $\lambda_{2} > 2\lambda_{1}$ separately. All the normal forms are discussed in Section \ref{Normal_forms_section}. For now, keep in mind that for all these normal forms the local stable and unstable invariant manifolds of $O$ are straightened, i.e. $W^{s}_{\mathrm{loc}}(O) = \{v_{1} = v_{2} = 0\}$ and $W^{u}_{\mathrm{loc}}(O) = \{u_{1} = u_{2} = 0\}$. Moreover, all the normal forms remain invariant with respect to symmetry (\ref{eq425}). In addition, reduction to the normal form preserves the form (\ref{First_integral}) of the first integral.

\begin{figure}[h]
\centering
\includegraphics[scale=0.09]{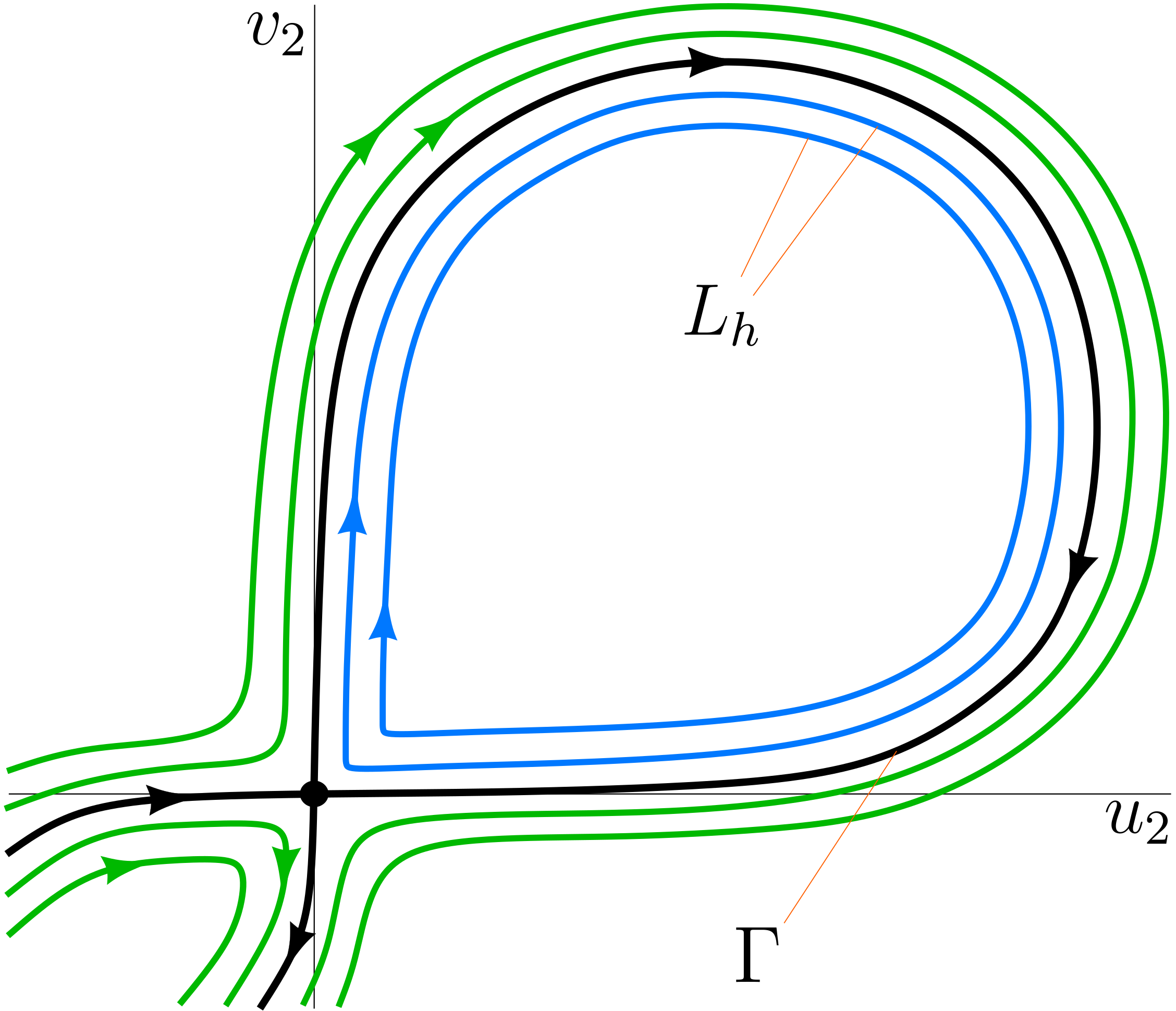}
\caption{\small The homoclinic orbit $\Gamma$ lies in the invariant $(u_{2}, v_{2})$-plane and the level set $\{H=0\}$. There exists a family of periodic orbits $L_h$ (shown by blue color) in the invariant plane close to $\Gamma$ in the region surrounded by the closed curve $\Gamma \cup \{O\}$. The periodic orbit $L_h$ lies in the level set $\{H=h\}$ for $-h_{0} \leq h < 0$, where $h_{0} > 0$ is sufficiently small.}
\label{Periodic_orbit_single_homoclinic_case}
\end{figure}

Restricting our system to the invariant plane $\{u_{1} = v_{1} = 0\}$, we obtain a planar system with the homoclinic $\Gamma$ to the equilibrium $(u_{2}, v_{2}) = (0,0)$. This system (at least in a small open neighborhood of $\Gamma$ in the plane) possesses a first integral (obtained from the restriction of $H$ to the plane). Therefore, there exists a continuum of periodic orbits $L_{h}$ inside the region surrounded by the closed loop $\Gamma\cup \{(0,0)\}$ in the plane, and accumulated to the loop (see Figure \ref{Periodic_orbit_single_homoclinic_case}). Each periodic loop $L_h$, where $-h_{0} \leq h < 0$ for some sufficiently small $h_{0} >0$, lies at the intersection of the level set $\{H = h\}$ and the invariant plane.

\begin{figure}
\centering
\includegraphics[scale=0.05]{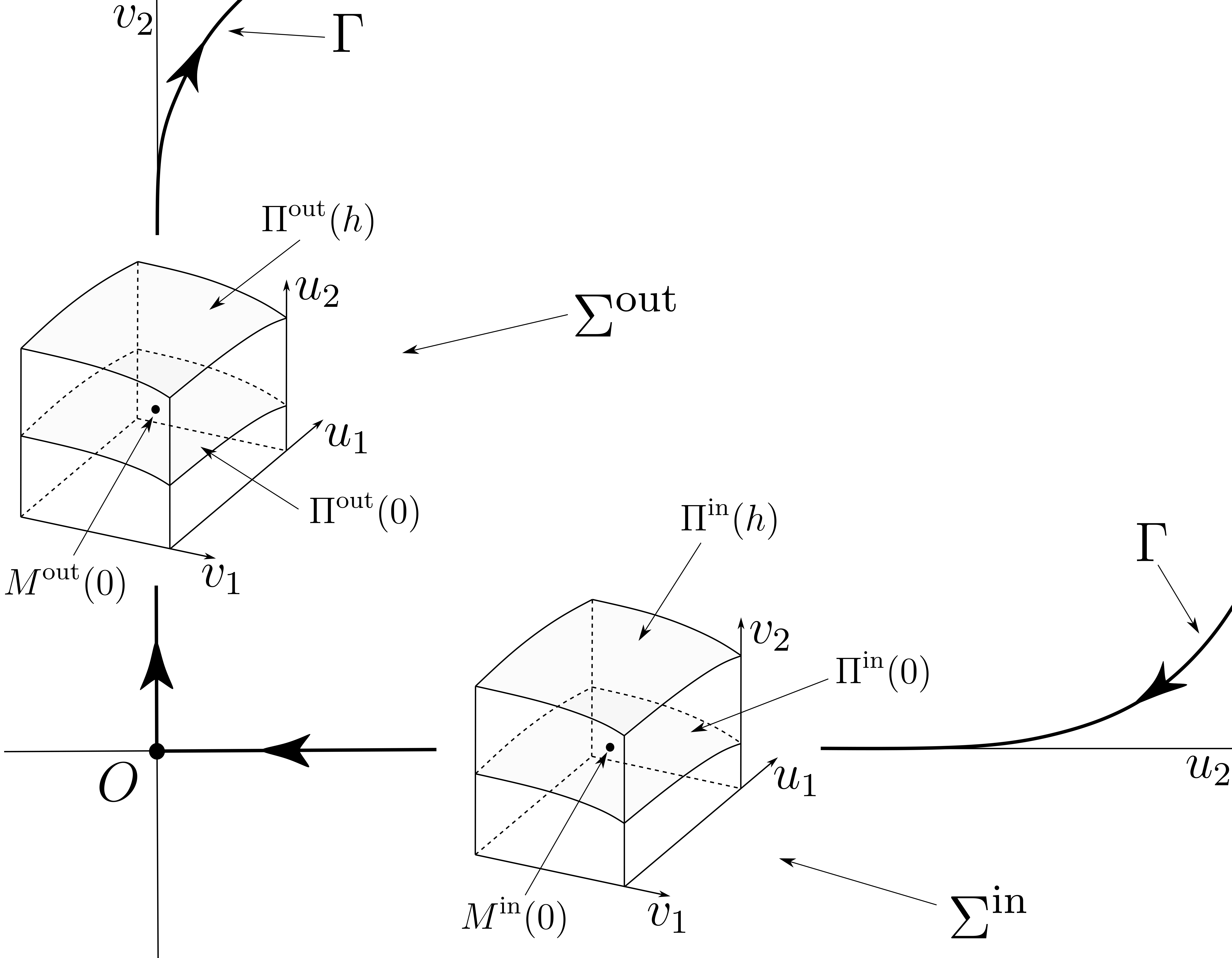}
\caption{\small The three-dimensional cross-sections $\Sigma^{\mathrm{in}}$ and $\Sigma^{\mathrm{out}}$ to the homoclinic orbit $\Gamma$ are shown. Each cross-section is foliated by two-dimensional sections $\Pi^{\mathrm{in/out}}(h)$.}
\label{Fig_three_dim_cross_sections}
\end{figure}

Consider two small 3-dimensional cross-sections $\Sigma^{\mathrm{in}} = \{u_2 = \delta\}$ and $\Sigma^{\mathrm{out}} = \{v_2 = \delta\}$ to the loop $\Gamma$, where $\delta > 0$ is sufficiently small (see Figure \ref{Fig_three_dim_cross_sections}). Take a sufficiently small $h_{0} > 0$, and let $\Pi^{\mathrm{in}}(h)$ and $\Pi^{\mathrm{out}}(h)$ be the restrictions of $\Sigma^{\mathrm{in}}$ and $\Sigma^{\mathrm{out}}$ to the level set $\{H = h\}$, respectively, i.e. \(\Pi^{\mathrm{in}}(h) = \lbrace u_{2} = \delta\rbrace \cap \lbrace H = h\rbrace\) and \(\Pi^{\mathrm{out}}(h) = \lbrace v_{2} = \delta\rbrace \cap \lbrace H = h\rbrace\). We can choose $(u_{1}, v_{1})$-coordinates on $\Pi^{\mathrm{in}}(h)$ and $\Pi^{\mathrm{out}}(h)$ meaning that $(u_{2}, v_{2})$ is uniquely determined by $(u_{1}, v_{1})$ (see Lemma \ref{Coordinates_Section_Lemma}).

The homoclinic orbit $\Gamma$ lies in the zero-level set of $H$. Thus, it intersects $\Sigma^{\mathrm{in}}$ at a point $M^{\mathrm{in}}(0)\in \Pi^{\mathrm{in}}(0)$, and $\Sigma^{\mathrm{out}}$ at a point $M^{\mathrm{out}}(0) \in \Pi^{\mathrm{out}}(0)$. As $\delta$ is sufficiently small, we have $M^{\mathrm{in}}(0) = (0,0)$ and $M^{\mathrm{out}}(0) = (0,0)$ (in the full $(u_{1}, u_{2}, v_{1}, v_{2})$-coordinates, we have that $M^{\mathrm{in}}(0) = (0, \delta, 0, 0)$ and $M^{\mathrm{out}}(0) = (0, 0, 0, \delta)$). Starting from $\Pi^{\mathrm{out}}(0)$, the forward orbit of any point sufficiently close to $M^{\mathrm{out}}(0)$ moves along $\Gamma$ and reaches $\Pi^{\mathrm{in}}(0)$ after a finite time. This defines a diffeomorphism $T_{0}^{\mathrm{glo}}$ from a small neighborhood of $M^{\mathrm{out}}(0)$ in $\Pi^{\mathrm{out}}(0)$ to a neighborhood of $M^{\mathrm{in}}(0)$ in $\Pi^{\mathrm{in}}(0)$. Taylor expanding $T^{\mathrm{glo}}_0$ at $M^{\mathrm{out}}(0)$, we have
\begin{equation}\label{eq63000}
T_{0}^{\text{glo}}\left(\begin{matrix}
u_1 \\ v_1
\end{matrix}\right) = \left(\begin{matrix}
a & b\\ c & d
\end{matrix}\right) \left(\begin{matrix}
u_1 \\ v_1
\end{matrix}\right) + \left(\begin{matrix}
o\left(\lvert u_{1}\rvert, \lvert v_{1}\rvert\right)\\
o\left(\lvert u_{1}\rvert, \lvert v_{1}\rvert\right)
\end{matrix}\right)
\end{equation}
for some \(a, b, c, d \in \mathbb{R}\). This map is the global piece (restricted to $\{H=0\}$) of the Poincar\'e map along $\Gamma$, which is discussed in more detail in Section \ref{Sec_Proof_Thm_Single_Loop}. Note that the transversality condition in Assumption \ref{assumption50} is equivalent to \(d\neq 0\). This is because the local unstable manifold of $O$ intersects $\Pi^{\mathrm{out}}(0)$ at $v_1$-axis and its image under $T^{\mathrm{glo}}$ must intersect the $u_1$-axis (which is the intersection of $\Pi^{\mathrm{in}}(0)$ and the local stable manifold of $O$) transversely. The following theorem is proved in \cite{Bakrani2022JDE}.

\begin{theorem}\label{Thm_JDE_1}(\cite{Bakrani2022JDE})
Assume $\lambda_{2}\neq 2\lambda_{1}$. For a given open neighborhood $U_0$ of $\Gamma \cup \{O\}$ in the level set $\{H=0\}$, let $W^{s}_{\mathrm{loc}}(\Gamma)$ (resp. $W^{u}_{\mathrm{loc}}(\Gamma)$) be the union of $\Gamma$ and the set of all points in $U_0$ whose forward (resp. backward) orbits lie entirely in $U_0$, and their $\omega$-limit (resp. $\alpha$-limit) sets coincide with $\Gamma \cup \{O\}$. Moreover, let $W^{s}_{U_0}(O)$ (resp. $W^{u}_{U_0}(O)$) be the set of the points on the global stable (resp. unstable) manifold of $O$ whose forward (resp. backward) orbits lie entirely in $U_0$. Assume \(b\), \(c\) and \(d\) in (\ref{eq63000}) are nonzero. Then, there exists a sufficiently small neighborhood $U_0$ such that the forward (resp. backward) orbit of a point in $U_0$ remains in $U_0$ if and only if it belongs to $W^{s}_{\mathrm{loc}}(\Gamma) \cup W^{s}_{U_0}(O)$ (resp. $W^{u}_{\mathrm{loc}}(\Gamma) \cup W^{u}_{U_0}(O)$). Moreover, the following hold.
\begin{enumerate}[(i)]
\item Suppose \(\lambda_{2} < 2\lambda_{1}\). Then, \(W^{s}_{\text{loc}}\left(\Gamma\right) = W^{u}_{\text{loc}}\left(\Gamma\right) = \Gamma\).
\item Suppose \(2\lambda_{1} < \lambda_{2}\). If \(cd > 0\) (resp. $bd < 0$), then \(W^{s}_{\text{loc}}(\Gamma)= \Gamma\) (resp. \(W^{u}_{\text{loc}}(\Gamma)= \Gamma\)). If \(cd<0\) (resp. $bd > 0$), then \(W^{s}_{\text{loc}}\left(\Gamma\right)\) (resp. \(W^{u}_{\text{loc}}\left(\Gamma\right)\)) is a \(\mathcal{C}^1\)-smooth 2-dimensional invariant manifold which is tangent to \(W_{\text{glo}}^{s}\left(O\right)\) (resp. \(W_{\text{glo}}^{u}\left(O\right)\)) at every point of \(\Gamma\).
\end{enumerate}
\end{theorem}

This theorem describes the dynamics near $\Gamma$ inside the level set $\{H = 0\}$. Roughly speaking, it states that an orbit in the level set $\{H = 0\}$ that starts close to $\Gamma$ and as $t\rightarrow\infty$ (resp. $t\rightarrow -\infty$), never leaves a small neighborhood of $\Gamma$ in $\{H = 0\}$ must belong to the stable (resp. unstable) set of $\Gamma$ or the stable (resp. unstable) manifold of $O$. In this paper, we are interested in describing the dynamics near $\Gamma$ in the nonzero level sets. Thus, together with the results of \cite{Bakrani2022JDE}, our theorems describe the dynamics in an open neighborhood of $\Gamma$ in the four-dimensional phase space $\mathbb{R}^{4}$.

For a given open neighborhood $U$ of $\Gamma \cup \{O\}$ in $\mathbb{R}^4$, let $U_{h} := U \cap \{H=h\}$. Obviously, $U_{h_{1}}\cap U_{h_2} = \emptyset$ when $h_{1} \neq h_{2}$. Recall the periodic orbit $L_{h}$ (see Figure \ref{Periodic_orbit_single_homoclinic_case}). Our first main result is the following.

\begin{theoremA}\label{Thm_Single_Loop}
Suppose $a$, $b$, $c$ and $d$ in (\ref{eq63000}) are nonzero, and Assumptions \ref{assumption30} and \ref{assumption50} hold. There exists a sufficiently small neighborhood $U$ of $\Gamma \cup \{O\}$, and a sufficiently small constant $h_0 > 0$ such that the following hold.
\begin{enumerate}[(i)]
\item If $-h_0 \leq h < 0$, the periodic orbit $L_h$ that lies in the invariant plane $\{u_1 = v_1 = 0\}$ is the only periodic orbit in $U_h$. Restricting system (\ref{eq74wt366cuw6gt5uw6v01092}) to the level set $\{H=h\}$, the periodic orbit $L_h$ is saddle possessing two-dimensional stable and unstable invariant manifolds $W^{s}(L_h)$ and $W^{u}(L_h)$, respectively. The forward (resp. backward) orbit of any point in $U_h$ off the manifold $W^{s}(L_h)$ (resp. $W^{u}(L_h)$) leaves $U_{h}$.
\item If $h=0$ and $\lambda_{2} \neq 2\lambda_{1}$, then the dynamics in $U_{0}$ is the one that is described by Theorem \ref{Thm_JDE_1}.
\item If $0 < h\leq h_0$, the forward and backward orbits of all the points in $U_h$ leave $U_h$.
\end{enumerate}
\end{theoremA}

The periodic orbit $L_h$ is saddle when the system is restricted to the level set $\{H=h\}$. Considering the whole 4-dimensional phase space, $L_h$ possesses a 2-dimensional local center manifold, which is indeed a neighborhood of it in the invariant $(u_{2}, v_{2})$-plane.

\begin{figure}[h]
\centering
\includegraphics[scale=.18]{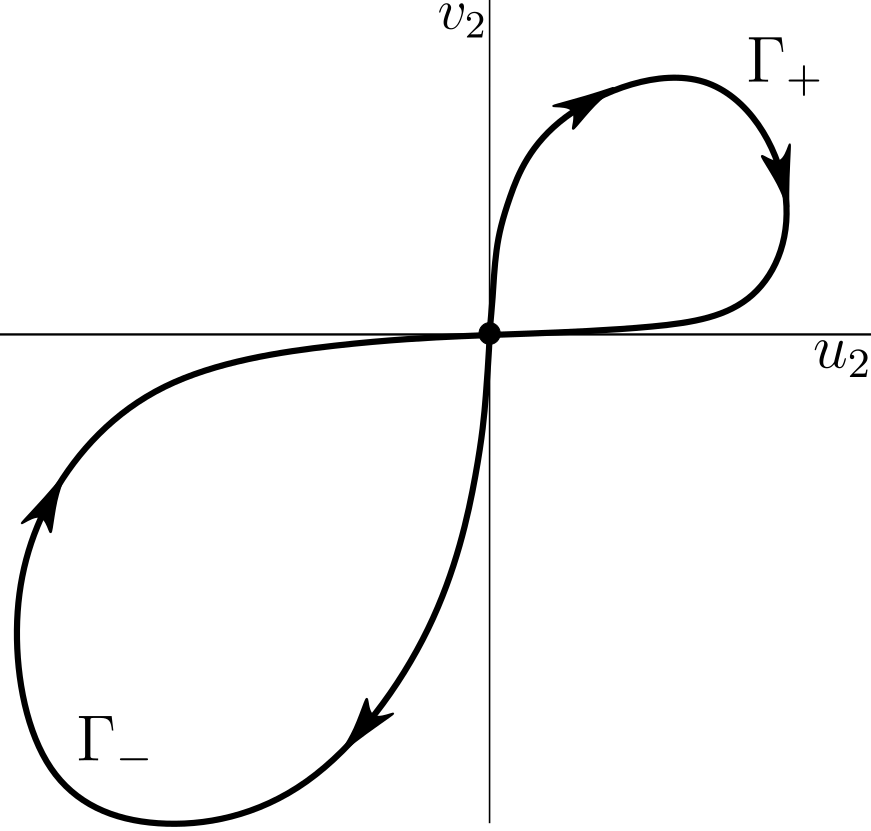}
\caption{\small A pair of transverse homoclinic loops in the invariant plane \(\lbrace u_{1}=v_{1}=0 \rbrace\).}
\label{Figureyuybbilzqbe564}
\end{figure}

So far, we have considered the case of a single homoclinic orbit. Let us now consider a second scenario in which we have two homoclinic orbits in the invariant $(u_{2},v_{2})$-plane.
\begin{myassumption}\label{assumption60}
There exist two homoclinic orbits \(\Gamma_{+}\) and \(\Gamma_{-}\) of the transverse intersection of \(W^{s}\left(O\right)\) and \(W^{u}\left(O\right)\) in the invariant plane \(\lbrace u_{1}=v_{1}=0 \rbrace\). The homoclinic orbit $\Gamma_{+}$ (resp. $\Gamma_{-}$) leaves and enters $O$ along the positive (resp. negative) sides of $v_2$ and $u_{2}$ axes, respectively (see Figure \ref{Figureyuybbilzqbe564}).
\end{myassumption}

We do not require $\Gamma_{+}$ and $\Gamma_{-}$ to be symmetric. The existence of a pair of homoclinic orbits (homoclinic figure-eight) is indeed a generic scenario when the level set \(\lbrace H = 0\rbrace\) is compact.

\begin{figure}[h]
\centering
\includegraphics[scale=.05]{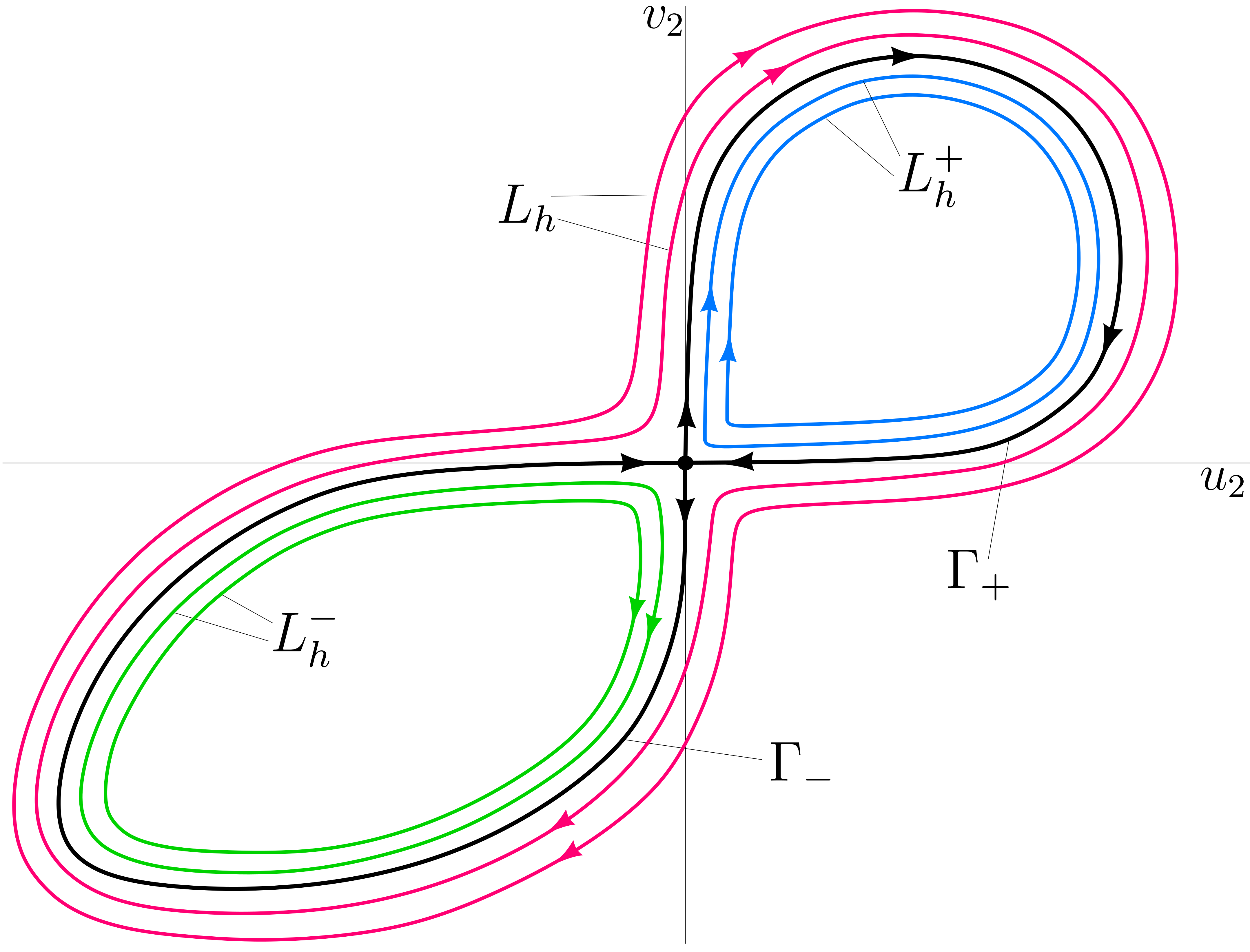}
\caption{\small The homoclinic orbits $\Gamma_{+}$ and $\Gamma_{-}$ lie in the invariant $(u_{2}, v_{2})$-plane. A continuum of periodic orbits $L_{h}^{+}$ (resp. $L_{h}^{-}$) shown by blue (resp. green) lies in the invariant plane inside the region surrounded by $\Gamma_{+} \cup \{O\}$ (resp. $\Gamma_{-} \cup \{O\}$). Moreover, a continuum of periodic orbits lies in the invariant plane (shown by red) outside the regions surrounded by $\Gamma_{+} \cup \{O\} \cup \Gamma_{-}$. In this figure, the homoclinic orbits are straightened near the equilibrium $O$. This does not need to be the case in general, however, this is the case when the system is brought into the normal forms (see Section \ref{Normal_forms_section}).}
\label{Periodic_orbits_double_homoclinic_case}
\end{figure}

Analogous to the case of a single homoclinic, there exist three families of periodic orbits $\{L_{h}\}_{0<h < h_{0}}$, $\{L^{+}_{h}\}_{-h_{0} <h < 0}$ and $\{L^{-}_{h}\}_{-h_{0} <h < 0}$ near $\Gamma_{+} \cup \{O\} \cup \Gamma_{-}$ lying in the invariant plane $\{u_{1} = v_{1} = 0\}$. The periodic orbits $L^{+}_{h}$ (resp. $L^{-}_{h}$) lie inside the region surrounded by $\Gamma_{+} \cup \{O\}$ (resp. $\Gamma_{-} \cup \{O\}$), and the periodic orbits $L_{h}$ lie outside of the region surrounded by $\Gamma_{+} \cup \{O\} \cup \Gamma_{-}$. Each of the orbits $L^{+}_h$, $L^{-}_h$, and $L_h$ lies in the level set $\{H=h\}$ (see Figure \ref{Periodic_orbits_double_homoclinic_case}).

\begin{figure}
\centering
\includegraphics[scale=0.035]{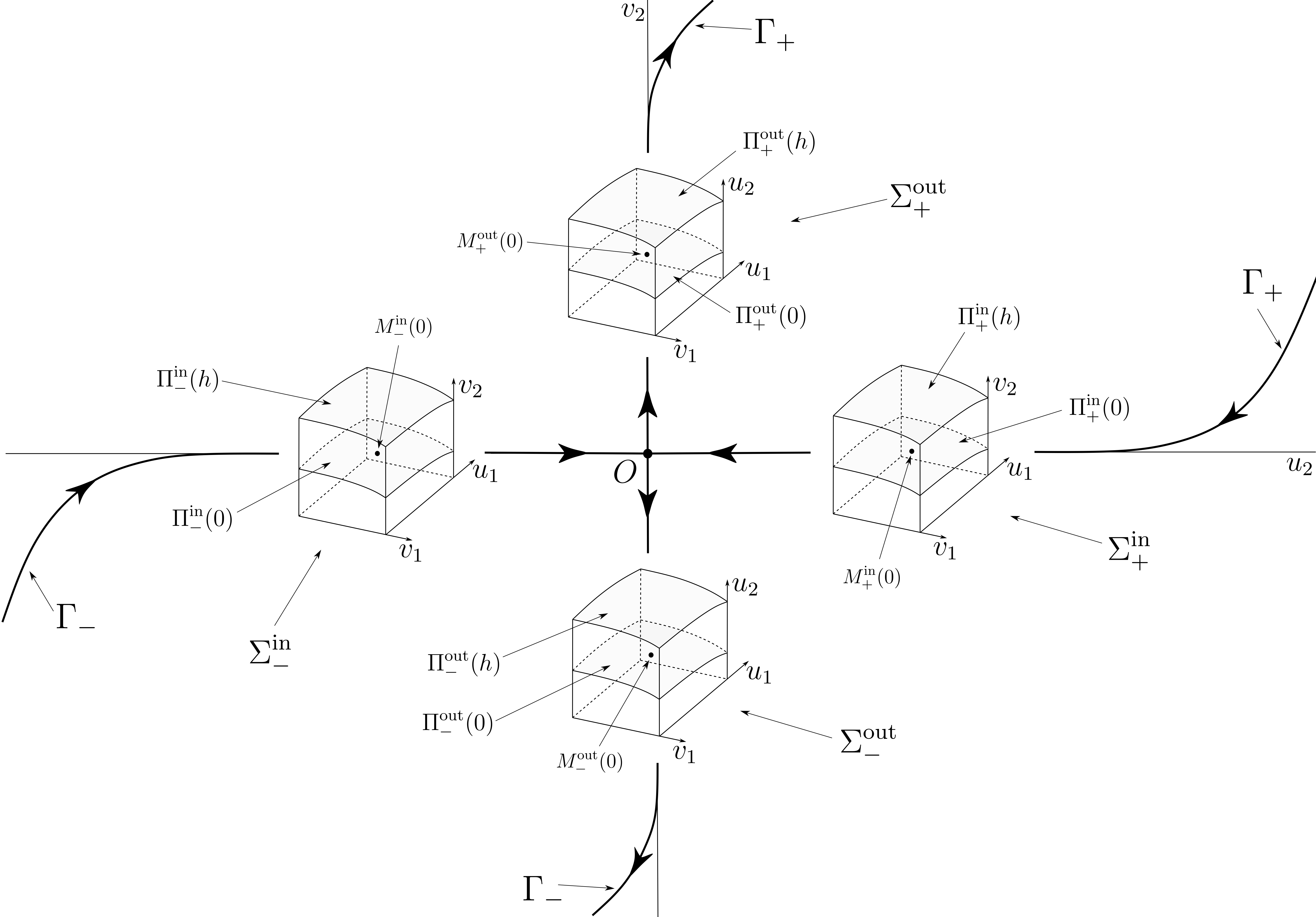}
\caption{\small The three-dimensional cross-sections $\Sigma_{\sigma}^{\mathrm{in}}$ and $\Sigma_{\sigma}^{\mathrm{out}}$ to the homoclinic orbit $\Gamma_{\sigma}$ for $\sigma = \pm$ are shown. Each cross-section is foliated by two-dimensional sections $\Pi_{\sigma}^{\mathrm{in/out}}(h)$.}
\label{Fig_three_dim_cross_sections_double_homoclinics}
\end{figure}

We aim to describe the dynamics near $\Gamma_{+}\cup \{O\} \cup \Gamma_{-}$ in the level set $\{H = h\}$, where $-h_0 \leq h \leq h_0$ for some small $h_0 > 0$. To this end, consider system (\ref{eq74wt366cuw6gt5uw6v01092}) again. Similar to the case of a single homoclinic $\Gamma$, suppose that the system near the origin is brought to the normal form (see Section \ref{Normal_forms_section}). We consider two small cross-sections $\Sigma^{\mathrm{in}}_{+} = \{u_{2} = \delta\}$ and $\Sigma^{\mathrm{out}}_{+} = \{v_{2} = \delta\}$ to the loop $\Gamma_{+}$, and two cross-sections $\Sigma^{\mathrm{in}}_{-} = \{u_{2} = -\delta\}$ and $\Sigma^{\mathrm{out}}_{-} = \{v_{2} = -\delta\}$ to the loop $\Gamma_{-}$. We define $\Pi^{\mathrm{in}}_{\pm}(h) = \Sigma^{\mathrm{in}}_{\pm} \cap \{H = h\}$ and $\Pi^{\mathrm{out}}_{\pm}(h) = \Sigma^{\mathrm{out}}_{\pm} \cap \{H = h\}$ (see Figure \ref{Fig_three_dim_cross_sections_double_homoclinics}). Note that at the intersection points of the homoclinics $\Gamma_{+}$ and $\Gamma_{-}$ with the cross-sections, we have $u_1 = v_1 = 0$. We can again choose $(u_1, v_1)$-coordinates on each of $\Pi_{\pm}^{\mathrm{in}}(h)$ and $\Pi_{\pm}^{\mathrm{out}}(h)$. Analogous to $a$, $b$, $c$ and $d$ in (\ref{eq63000}), let $a_\sigma$, $b_\sigma$, $c_\sigma$ and $d_\sigma$ for $\sigma = \pm$ be the associated coefficients of the global map along $\Gamma_\sigma$. The following theorem from \cite{Bakrani2022JDE} holds.

\begin{theorem}\label{Thm_JDE_2}(\cite{Bakrani2022JDE})
Assume $\lambda_{2}\neq 2\lambda_{1}$. Consider $W^{s}_{\mathrm{loc}}(\Gamma_{+})$, $W^{s}_{\mathrm{loc}}(\Gamma_{-})$, $W^{u}_{\mathrm{loc}}(\Gamma_{+})$ and $W^{u}_{\mathrm{loc}}(\Gamma_{-})$ given by Theorem \ref{Thm_JDE_1}. We also define $W^{s}_{\mathrm{loc}}(\Gamma_{+}\cup\Gamma_{-})$ (resp. $W^{u}_{\mathrm{loc}}(\Gamma_{+}\cup\Gamma_{-})$) to be the union of $\Gamma_{+}\cup\Gamma_{-}$ and the set of all points in an open neighborhood $V_0$ of $\Gamma_{+} \cup \{O\}\cup\Gamma_{-}$ in the level set $\{H=0\}$ whose forward (resp. backward) orbits lie entirely in $V_0$, and their $\omega$-limit (resp. $\alpha$-limit) sets coincide with $\Gamma_{+} \cup \{O\}\cup \Gamma_{-}$. Define $W^{s}_{V_0}(O)$ and $W^{u}_{V_0}(O)$ analogous to the ones in Theorem \ref{Thm_JDE_1}. Assume \(b_\sigma\), \(c_\sigma\) and \(d_\sigma\) for $\sigma = \pm$ are nonzero. Then, there exists a sufficiently small neighborhood $V_0$ such that the forward (resp. backward) orbit of a point in $V_0$ remains in $V_0$ if and only if it belongs to $W^{s}_{\mathrm{loc}}(\Gamma_{+}) \cup W^{s}_{\mathrm{loc}}(\Gamma_{-}) \cup W^{s}_{\mathrm{loc}}(\Gamma_{+}\cup \Gamma_{-}) \cup W^{s}_{V_0}(O)$ (resp. $W^{u}_{\mathrm{loc}}(\Gamma_{+}) \cup W^{u}_{\mathrm{loc}}(\Gamma_{-}) \cup W^{u}_{\mathrm{loc}}(\Gamma_{+}\cup \Gamma_{-}) \cup W^{u}_{V_0}(O)$). Moreover, the following hold.
\begin{enumerate}[(i)]
\item Suppose \(\lambda_{2} < 2\lambda_{1}\). Then, \(W^{s}_{\mathrm{loc}}\left(\Gamma_{+}\cup\Gamma_{-}\right) = \Gamma_{+}\cup\Gamma_{-} = W^{u}_{\mathrm{loc}}\left(\Gamma_{+}\cup\Gamma_{-}\right)\).
\item Suppose \(2\lambda_{1} < \lambda_{2}\). If \(c_{+}d_{+} > 0\) and \(c_{-}d_{-} > 0\) (resp. $b_{+}d_{+} < 0$ and $b_{-}d_{-} < 0$), then \(W^{s}_{\mathrm{loc}}(\Gamma_{+} \cup \Gamma_{-})\) (resp. \(W^{u}_{\mathrm{loc}}(\Gamma_{+} \cup \Gamma_{-})\)) is a \(\mathcal{C}^1\)-smooth 2-dimensional invariant manifold which is tangent to \(W_{\mathrm{glo}}^{s}\left(O\right)\) (resp. \(W_{\mathrm{glo}}^{u}\left(O\right)\)) at every point of \(\Gamma_{+} \cup \Gamma_{-}\). Otherwise, $W^{s}_{\mathrm{loc}}(\Gamma_{+} \cup \Gamma_{-}) = \Gamma_{+} \cup \Gamma_{-}$ (resp. $W^{u}_{\mathrm{loc}}(\Gamma_{+} \cup \Gamma_{-}) = \Gamma_{+} \cup \Gamma_{-}$).
\end{enumerate}
\end{theorem}

This theorem describes the dynamics near $\Gamma_{+} \cup \{O\} \cup\Gamma_{-}$ in the level set $\{H=0\}$. The following theorem, which is our second main result, describes the dynamics near $\Gamma_{+} \cup \{O\} \cup\Gamma_{-}$ in nonzero level sets. To state the theorem, consider an  open neighborhood $V$ of $\Gamma_{+} \cup \{O\} \cup\Gamma_{-}$ in $\mathbb{R}^4$, and let $V_{h} := V \cap \{H=h\}$. Obviously, $V_{h_{1}}\cap V_{h_2} = \emptyset$ for distinct $h_{1}, h_{2}\in \mathbb{R}$. Recall the periodic orbits $L_{h}$, $L^{+}_{h}$ and $L^{-}_{h}$ (see Figure \ref{Periodic_orbits_double_homoclinic_case}). Then, the following holds.

\begin{theoremB}\label{Thm_Double_Loops}
Suppose $a_\sigma$, $b_\sigma$, $c_\sigma$ and $d_\sigma$ for $\sigma = \pm$ are nonzero, and Assumptions \ref{assumption30} and \ref{assumption60} hold. There exists a sufficiently small neighborhood $V$ of $\Gamma_{+} \cup \{O\} \cup \Gamma_{-}$, and a sufficiently small constant $h_0 > 0$ such that the following hold.
\begin{enumerate}[(i)]
\item \label{ItemThmDouble0001} If $-h_{0} \leq h < 0$, the periodic orbits $L_h^{+}$ and $L_h^{-}$ that lie in the invariant plane $\{u_1 = v_1 = 0\}$ are the only periodic orbits in $V_h$. Restricting system (\ref{eq74wt366cuw6gt5uw6v01092}) to the level set $\{H=h\}$, the periodic orbits $L_h^{+}$ and $L_h^{-}$ are saddle possessing two-dimensional stable and unstable invariant manifolds $W^{s}(L_h^{\pm})$ and $W^{u}(L_h^{\pm})$, respectively. The forward (resp. backward) orbit of any point in $V_h$ off the set $W^{s}(L_h^{+}) \cup W^{s}(L_h^{-})$ (resp. $W^{u}(L_h^{+}) \cup W^{u}(L_h^{-})$) leaves $V_{h}$.
\item If $h=0$ and $\lambda_{2} \neq 2\lambda_{1}$, then the dynamics in $V_{0}$ is the one that is described by Theorem \ref{Thm_JDE_2}.
\item \label{ItemThmDouble0002} If $0 < h \leq h_{0}$, the periodic orbit $L_h$ that lies in the invariant plane $\{u_1 = v_1 = 0\}$ is the only periodic orbit in $V_h$. Restricting system (\ref{eq74wt366cuw6gt5uw6v01092}) to the level set $\{H=h\}$, the periodic orbit $L_h$ is saddle possessing two-dimensional stable and unstable invariant manifolds $W^{s}(L_h)$ and $W^{u}(L_h)$, respectively. The forward (resp. backward) orbit of any point in $V_h$ off the manifold $W^{s}(L_h)$ (resp. $W^{u}(L_h)$) leaves $V_{h}$.
\end{enumerate}
\end{theoremB}

All the periodic orbits $L_h$, $L^{+}_{h}$ and $L^{-}_h$ are saddle when the system is restricted to $\{H=h\}$. Considering the whole 4-dimensional phase space, all of them possess 2-dimensional local center manifolds (each center manifold is a neighborhood of the periodic orbit in the invariant $(u_{2}, v_{2})$-plane).

The proofs of our main results are based on the study of the Poincar\'e maps along the homoclinic and periodic orbits. To this end, we first reduce the system near the origin to a normal form. This is done in Section \ref{Normal_forms_section}. Then, in Section \ref{Section_trajectories_near_O}, we approximate the trajectories near the origin. This allows us to approximate the Poincar\'e maps and prove Theorem \ref{Thm_Single_Loop} in Section \ref{Sec_Proof_Thm_Single_Loop}, and Theorem \ref{Thm_Double_Loops} in Section \ref{Sec_Proof_Thm_Double_Loops}. We leave the proof of some of the technical statements for Section \ref{Technical_Section}.

\section{Normal forms}\label{Normal_forms_section}

To study the Poincar\'e maps along the homoclinic and periodic orbits, we first bring system (\ref{eq300}) near the origin to a normal form. This is done in this section for three cases \(\lambda_{1} = \lambda_{2}\), \(\lambda_{1} < \lambda_{2} \leq 2\lambda_{1}\) and \(2\lambda_{1} < \lambda_{2}\) seperately. The normal forms that are used here are the same as the ones in \cite{Bakrani2022JDE} (see Section 3.2.2 of \cite{BakraniPhDthesis} for the proofs). Despite \cite{Bakrani2022JDE} that the case $\lambda_{2} = 2\lambda_{1}$ was excluded, we study this resonant case in this paper. The normal form that we use for this case is the same as the one for the case $\lambda_{1} < \lambda_{2} < 2\lambda_{1}$ (see Lemma \ref{Nftheorem} below).

\begin{mylem}[Normal form for the case $\lambda_{1} = \lambda_{2}$]\label{Nftheorem0}
Consider system (\ref{eq300}) and first integral (\ref{eq400}) for $\lambda := \lambda_{1} = \lambda_{2}$. There exists a \({\mathcal{C}}^{\infty}\)-smooth change of coordinates which brings system (\ref{eq300}) to the form
\begin{equation}\label{eq74wt366cuw6gt5uw6v01092}
\begin{aligned}
\dot{u}_{1} &= -\lambda u_{1} + f_{11}(u_{1}, u_{2}, v_{1}, v_{2}) u_{1} + f_{12}(u_{1}, u_{2}, v_{1}, v_{2}) u_{2},\\
\dot{u}_{2} &= -\lambda u_{2} + f_{21}(u_{1}, u_{2}, v_{1}, v_{2}) u_{1} + f_{22}(u_{1}, u_{2}, v_{1}, v_{2}) u_{2},\\
\dot{v}_{1} &= +\lambda v_{1} + g_{11}(u_{1}, u_{2}, v_{1}, v_{2}) v_{1} + g_{12}(u_{1}, u_{2}, v_{1}, v_{2}) v_{2},\\
\dot{v}_{2} &= +\lambda v_{2} + g_{21}(u_{1}, u_{2}, v_{1}, v_{2}) v_{1} + g_{22}(u_{1}, u_{2}, v_{1}, v_{2}) v_{2},
\end{aligned}
\end{equation}
where the functions \(f_{ij}\), \(g_{ij}\) are \({\mathcal{C}}^{\infty}\)-smooth and vanish at the origin, i.e. $f_{ij}\left(0,0,0,0\right) = g_{ij}\left(0,0,0,0\right) = 0$, and transforms first integral (\ref{eq400}) to $H = \lambda u_{1}v_{1} - \lambda u_{2}v_{2}$. By a rescaling, we set
\begin{equation}\label{eq98n9uo8wbr8bau4b819ub92q1u8}
H = u_{1}v_{1} - u_{2}v_{2}.
\end{equation}
Moreover, system (\ref{eq74wt366cuw6gt5uw6v01092}) remains invariant with respect to symmetry (\ref{eq425}). In particular,
\begin{equation}\label{eq24500}
f_{12}(0, u_{2}, 0, v_{2}) \equiv 0,\quad g_{12}(0, u_{2}, 0, v_{2}) \equiv 0.
\end{equation}
\end{mylem}

\begin{mylem}[Normal form for the cases $\lambda_{1} < \lambda_{2} < 2\lambda_{1}$ and $\lambda_{2} = 2\lambda_{1}$]\label{Nftheorem}
Consider system (\ref{eq300}) and first integral (\ref{eq400}), and assume \(\lambda_{1} <\lambda_{2}\). There exists a \({\mathcal{C}}^{\infty}\)-smooth change of coordinates which brings system (\ref{eq300}) to the form
\begin{equation}\label{eq19000}
\begin{aligned}
\dot{u}_{1} &= -\lambda_{1}u_{1} + f_{11}\left(u_{1}, v\right) u_{1} + f_{12}\left(u_{1}, u_{2}, v\right) u_{2},\\
\dot{u}_{2} &= -\lambda_{2}u_{2} + f_{21}\left(u_{1}, v\right)u_{1} + f_{22}\left(u_{1}, u_{2}, v\right) u_{2},\\
\dot{v}_{1} &= +\lambda_{1}v_{1} + g_{11}\left(u, v_{1}\right) v_{1} + g_{12}\left(u, v_{1}, v_{2}\right) v_{2},\\
\dot{v}_{2} &= +\lambda_{2}v_{2} + g_{21}\left(u, v_{1}\right)v_{1} +  g_{22}\left(u, v_{1}, v_{2}\right) v_{2},
\end{aligned}
\end{equation}
where the functions \(f_{ij}\), \(g_{ij}\) are \({\mathcal{C}}^{\infty}\)-smooth and satisfy the identities
\begin{equation}\label{eq20000}
\begin{gathered}
f_{11}(0,v) \equiv 0,\quad f_{11}(u_{1},0) \equiv 0,\quad f_{12}(u,0) \equiv 0,\quad f_{21}(0,v) \equiv 0,\quad f_{22}(0, v) \equiv 0,\\
g_{11}(u,0) \equiv 0,\quad g_{11}(0,v_{1}) \equiv 0,\quad g_{12}(0,v) \equiv 0,\quad g_{21}(u,0) \equiv 0,\quad g_{22}\left(u,0\right) \equiv 0.
\end{gathered}
\end{equation}
This change of coordinates transforms first integral (\ref{eq400}) to
\begin{equation}\label{eq22000}
H = \lambda_{1}u_{1}v_{1} \left[1 + H_{1}\left(u, v\right)\right]-\lambda_{2}u_{2}v_{2} \left[1 + H_{2}\left(u, v\right)\right],
\end{equation}
where \(H_{1}\) and \(H_{2}\) are \(\mathcal{C}^{\infty}\) functions that vanish at \(O\). By a rescaling, we write (\ref{eq22000}) as
\begin{equation}\label{eq22000yuio}
H = \gamma u_{1}v_{1} \left[1 + o\left(1\right)\right] - u_{2}v_{2} \left[1 + o\left(1\right)\right], \quad \mathrm{where\,\,\,} \gamma =\lambda_{1} \lambda_{2}^{-1}.
\end{equation}
Moreover, normal form (\ref{eq19000}) and first integral (\ref{eq22000}) remain invariant with respect to symmetry (\ref{eq425}). In particular, (\ref{eq24500}) holds.
\end{mylem}

The statement of Lemma \ref{Nftheorem} holds for arbitrary \(\lambda_{1} <\lambda_{2}\). However, we will particularly use this normal form to analyze the local dynamics near \(O\) when \(\lambda_{1} <\lambda_{2} < 2\lambda_{1}\) or $\lambda_{2} = 2\lambda_{1}$. The normal form that is used for analyzing the case \(2\lambda_{1} < \lambda_{2}\) is given by the following:
\begin{mylem}[Normal form for the case $\lambda_{2} > 2\lambda_{1}$]\label{Nftheorem2}
Consider system (\ref{eq300}) and first integral (\ref{eq400}) and assume \(2\lambda_{1} <\lambda_{2}\). Let \(q\) be the largest integer such that \(q\lambda_{1} < \lambda_{2}\). There exists a \({\mathcal{C}}^{q}\)-smooth change of coordinates which brings system (\ref{eq300}) to the form
\begin{equation}\label{eq23000}
\begin{aligned}
\dot{u}_{1} &= -\lambda_{1}u_{1} + f_{11}\left(u_{1}, v\right) u_{1} + f_{12}\left(u_{1}, u_{2}, v\right) u_{2},\\
\dot{u}_{2} &= -\lambda_{2}u_{2} + f_{22}\left(u_{1}, u_{2}, v\right) u_{2},\\
\dot{v}_{1} &= +\lambda_{1}v_{1} + g_{11}\left(u, v_{1}\right) v_{1} + g_{12}\left(u, v_{1}, v_{2}\right) v_{2},\\
\dot{v}_{2} &= +\lambda_{2}v_{2} + g_{22}\left(u, v_{1}, v_{2}\right) v_{2},
\end{aligned}
\end{equation}
where \(f_{ij}\) and \(g_{ij}\) are \(\mathcal{C}^{q-1}\)-smooth and satisfy identities (\ref{eq20000}). This change of coordinates transforms first integral (\ref{eq400}) to
\begin{equation}\label{eqhhuhut5989b8b8y4qlz37y7t}
H = \lambda_{1}u_{1}v_{1} \left[1 + H_{1}\left(u, v\right)\right] -\lambda_{2} u_{2}v_{2} \left[1 + H_{2}\left(u, v\right)\right] + u_{2}v_{1}^{2} H_{3}\left(u, v\right) + v_{2}u_{1}^{2} H_{4}\left(u, v\right),
\end{equation}
where \(H\) is \(\mathcal{C}^{q}\), and \(H_{1}\), \(H_{2}\), \(H_{3}\) and \(H_{4}\) are some \(\mathcal{C}^{q-1}\), \(\mathcal{C}^{q}\), \(\mathcal{C}^{q-2}\) and \(\mathcal{C}^{q-2}\) functions, respectively, such that \(H_{1}(O) = H_{2}(O) = 0\). Moreover, system (\ref{eq23000}) and first integral (\ref{eqhhuhut5989b8b8y4qlz37y7t}) remain invariant with respect to symmetry (\ref{eq425}). In particular, (\ref{eq24500}) holds. By a rescaling, we write
\begin{equation}\label{eq12hniuwb7gydcw}
H = \gamma u_{1}v_{1} \left[1 + o\left(1\right)\right] - u_{2}v_{2} \left[1 + o\left(1\right)\right] + u_{2}v_{1}^{2} O\left(1\right) + v_{2}u_{1}^{2} O\left(1\right), \quad \mathrm{where\,\,\,} \gamma =\lambda_{1} \lambda_{2}^{-1}.
\end{equation}
\end{mylem}

\begin{myrem}
In all these three normal forms, the local stable and unstable as well as the local strong stable and strong unstable invariant manifolds of the equilibrium \(O\) are straightened, i.e. \(W^{s}_{\text{loc}} = \{v = 0\}\), \(W^{u}_{\text{loc}} = \{u = 0\}\), \(W^{ss}_{\text{loc}} = \{u_{1} = v_{1} = v_{2} = 0\}\) and \(W^{uu}_{\text{loc}} = \{u_{1} = u_{2} = v_{1} = 0\}\). In the case of normal form (\ref{eq23000}), the local extended stable and extended unstable invariant manifolds of \(O\) are straightened too, i.e. \(W^{sE}_{\text{loc}} = \{v_{2} = 0\}\) and \(W^{uE}_{\text{loc}} = \{u_{2} = 0\}\).
\end{myrem}

\section{Trajectories near the equilibrium \(O\)}\label{Section_trajectories_near_O}

To study the Poincar\'e map (particularly, the local map) along the homoclinics, we need to approximate the trajectories near the equilibrium $O$. Our method is based on the Shilnikov's boundary value problem (Shilnikov coordinates). Here, we first describe this method.

Consider a four-dimensional system
\begin{equation}\label{eqbwiebro4ibyw43yyyiw0290}
\begin{aligned}
\dot{u}_{1} &= -\lambda_{1} u_{1} + F_{1}\left(u_{1}, u_{2}, v_{1}, v_{2}\right),\qquad
&&\dot{v}_{1} = \lambda_{1} v_{1} + G_{1}\left(u_{1}, u_{2}, v_{1}, v_{2}\right),\\
\dot{u}_{2} &= -\lambda_{2} u_{2} + F_{2}\left(u_{1}, u_{2}, v_{1}, v_{2}\right),\qquad
&&\dot{v}_{2} = \lambda_{2} v_{2} + G_{2}\left(u_{1}, u_{2}, v_{1}, v_{2}\right),
\end{aligned}
\end{equation}
where $\lambda_{2} \geq \lambda_{1} >0$, and \(F_{1}\), \(F_{2}\), \(G_{1}\) and \(G_{2}\) and their first derivatives vanish at $O$. Take a $\delta>0$, and consider four real constants $u_{10}$, $u_{20}$, $v_{1\tau}$ and $v_{2\tau}$ such that $\max\{\lvert u_{10}\rvert, \lvert u_{20}\rvert, \lvert v_{1\tau}\rvert, \lvert v_{2\tau}\rvert\}\leq \delta$. If $\delta$ is sufficiently small, for any $\tau > 0$, there exists a unique solution
\begin{equation}\label{eq800050}
\left(u^{*}, v^{*}\right) = \left(u^{*}\left(t\right), v^{*}\left(t\right)\right) = \big(u^{*}\left(t, \tau, u_{10}, u_{20}, v_{1\tau}, v_{2\tau}\right), v^{*}\left(t, \tau, u_{10}, u_{20}, v_{1\tau}, v_{2\tau}\right)\big)
\end{equation}
of (\ref{eqbwiebro4ibyw43yyyiw0290}), where \(u^{*} = (u_{1}^{*},u_{2}^{*})\) and \(v^{*} = (v_{1}^{*},v_{2}^{*})\), that satisfies the boundary conditions
\begin{equation}\label{eq800020}
u^{*}_{1}\left(0\right) = u_{10}, \quad u^{*}_{2}\left(0\right) = u_{20}, \quad v^{*}_{1}\left(\tau\right) = v_{1\tau}, \quad v^{*}_{2}\left(\tau\right) = v_{2\tau}
\end{equation}
(see \cite{PoincareBirkhoff} and \cite{Dimabook}). The dependence of this solution on $(t, \tau, u_{10}, u_{20}, v_{1\tau}, v_{2\tau})$ is as smooth as system (\ref{eqbwiebro4ibyw43yyyiw0290}). Here, we describe a successive procedure to find this solution. To do this, observe that \(\left(u^{*}(t), v^{*}(t)\right)\) is a solution of (\ref{eqbwiebro4ibyw43yyyiw0290}) with boundary conditions (\ref{eq800020}) if and only if
\begin{equation}\label{eq800025}
\begin{aligned}
u_{1}^{*}(t) =& e^{-\lambda_{1} t} u_{10} + \int_{0}^{t} e^{\lambda_{1} (s-t)} F_{1}(x^{*}(s)) ds,\quad && v_{1}^{*}(t) = e^{-\lambda_{1} \left(\tau-t\right)} v_{1\tau} - \int_{t}^{\tau} e^{-\lambda_{1} (s-t)} G_{1}(x^{*}(s)) ds,\\
u_{2}^{*}(t) =& e^{-\lambda_{2} t} u_{20} + \int_{0}^{t} e^{\lambda_{2} (s-t)} F_{2}(x^{*}(s)) ds,\quad
&&v_{2}^{*}(t) = e^{-\lambda_{2} \left(\tau-t\right)} v_{2\tau} - \int_{t}^{\tau} e^{-\lambda_{2} (s-t)} G_{2}(x^{*}(s)) ds,
\end{aligned}
\end{equation}
where $x^{*}(s) := (u^{*}(s), v^{*}(s))$. Let $\Omega$ be the set of all the points $(u,v)$ such that $\max\{\lvert u_1\rvert, \lvert u_2\rvert, \lvert v_1\rvert, \lvert v_2\rvert\}\leq \delta$, for some $\delta >0$. Fix a $\tau > 0$, and let \(\mathcal{I}\) be the set of all the continuous functions
\begin{equation*}
\left[0, \tau\right] \ni t \mapsto \left(u_{1}\left(t\right), u_{2}\left(t\right), v_{1}\left(t\right), v_{2}\left(t\right)\right) \in \Omega.
\end{equation*}
Then, the integral operator $\mathfrak{T} : \mathcal{I}\rightarrow\mathcal{I}$ defined by
\begin{equation*}
\left(u_{1}\left(t\right), u_{2}\left(t\right), v_{1}\left(t\right), v_{2}\left(t\right)\right) \mapsto \left(\overline{u}_{1}\left(t\right), \overline{u}_{2}\left(t\right), \overline{v}_{1}\left(t\right), \overline{v}_{2}\left(t\right)\right),
\end{equation*}
where
\begin{equation*}
\begin{aligned}
\overline{u}_{1}\left(t\right) &= e^{-\lambda_{1} t} u_{10} + \int_{0}^{t} e^{\lambda_{1} (s-t)} F_{1}(x(s)) ds,\quad
&&\overline{v}_{1}\left(t\right) = e^{-\lambda_{1} \left(\tau -t\right)} v_{1\tau} - \int_{t}^{\tau} e^{-\lambda_{1} (s-t)} G_{1}(x(s)) ds,\\
\overline{u}_{2}\left(t\right) &= e^{-\lambda_{2} t} u_{20} + \int_{0}^{t} e^{\lambda_{2} (s-t)} F_{2}(x(s)) ds,\quad
&&\overline{v}_{2}\left(t\right) = e^{-\lambda_{2} \left(\tau -t\right)} v_{2\tau} - \int_{t}^{\tau} e^{-\lambda_{2} (s-t)} G_{2}(x(s)) ds,
\end{aligned}
\end{equation*}
for $x(s) = (u(s), v(s))$, is well-defined and a contraction on $\mathcal{I}$ \cite{Dimabook}. The set $\mathcal{I}$ endowed with $\mathcal{C}^{0}$-norm is a complete metric space. Therefore, thanks to the contraction mapping principle, $\mathfrak{T}$ has a unique fixed point, and iterating $\mathfrak{T}$ starting from $(u_{1}(t), u_{2}(t), v_{1}(t), v_{2}(t)) \equiv (0,0,0,0)$ converges to this fixed point. However, by (\ref{eq800025}), the fixed point of $\mathfrak{T}$ is nothing but solution (\ref{eq800050}).

Using this method, the trajectories near the origin of the normal forms of the cases $\lambda_{1} = \lambda_{2}$, $\lambda_{1}< \lambda_{2} < 2\lambda_{1}$ and $\lambda_{2} > 2\lambda_{1}$ are estimated in \cite{BakraniPhDthesis}. For the two cases $\lambda_{1}< \lambda_{2} < 2\lambda_{1}$ and $\lambda_{2} > 2\lambda_{1}$, we use the same estimates in this paper. For the case $\lambda_{1} = \lambda_{2}$, we find new sharper estimates below. We further find estimates for the resonant case $\lambda_{2} = 2\lambda_{1}$ below which was not studied in \cite{BakraniPhDthesis}.

\begin{mylem}[Case $\lambda_{1} = \lambda_{2}$]\label{flowlemmaresonant}
Let \(\lambda = \lambda_{1} = \lambda_{2}\). There exists \(M > 0\) such that for any sufficiently small \(\delta > 0\), and any \(u_{10}\), \(u_{20}\), \(v_{1\tau}\) and \(v_{2\tau}\), where \(\max \lbrace \lvert u_{10}\rvert, \lvert u_{20}\rvert, \lvert v_{1\tau}\rvert, \lvert v_{2\tau}\rvert\rbrace \leq \delta\), the solution \(\left(u\left(t\right), v\left(t\right)\right)\) of system (\ref{eq74wt366cuw6gt5uw6v01092}) that satisfies boundary conditions (\ref{eq800020}) can be written as
\begin{equation}\label{eq7bwwkauser6518092hz9ub}
\begin{aligned}
u_{1}(t) =& e^{-\lambda t}u_{10} + \xi_{1}\left(x\right),\qquad
&& u_{2}(t) = e^{-\lambda t}u_{20} + \xi_{2}\left(x\right),\\
v_{1}(t) =& e^{-\lambda\left(\tau - t\right)} v_{1\tau} + \zeta_{1}\left(x\right),\qquad
&& v_{2}(t) = e^{-\lambda\left(\tau-t\right)} v_{2\tau} +\zeta_{2}\left(x\right),
\end{aligned}
\end{equation}
where \(x = \left(t, \tau, u_{10}, u_{20}, v_{1\tau}, v_{2\tau}\right)\), \(t\in\left[0, \tau\right]\), and
\begin{equation}\label{eq879ygo7f6gvjhferg}
\begin{aligned}
\lvert \xi_{1}\rvert \leq & M \left(e^{-\lambda t}\lvert u_{10}\rvert\delta + e^{-\lambda\tau} \lvert v_{1\tau}\rvert \delta\right),\qquad &&
\lvert \xi_{2}\rvert \leq M e^{-\lambda t} \delta^{2},\\
\lvert \zeta_{1}\rvert \leq & M \left(e^{-\lambda (\tau -t)}\lvert v_{1\tau}\rvert \delta + e^{-\lambda\tau} \lvert u_{10}\rvert \delta\right),\qquad &&
\lvert \zeta_{2}\rvert \leq M e^{-\lambda (\tau - t)} \delta^{2}.
\end{aligned}
\end{equation}
\end{mylem}

\begin{myrem}
By (\ref{eq879ygo7f6gvjhferg}), the solution \(\left(u\left(t\right), v\left(t\right)\right)\) of system (\ref{eq74wt366cuw6gt5uw6v01092}) can be written as
\begin{equation}\label{eq65f3o8fgrcdc}
\begin{aligned}
u_{1}(t) &= e^{-\lambda t} u_{10}\left[1 + O\left(\delta\right)\right] + e^{-\lambda\tau} O\left(\lvert v_{1\tau}\rvert\delta\right),\qquad && u_{2}(t) = e^{-\lambda t} \left[u_{20} + O\left(\delta^{2}\right)\right],\\
v_{1}(t) &= e^{-\lambda (\tau -t)} v_{1\tau}\left[1 + O\left(\delta\right)\right] + e^{-\lambda\tau} O\left(\lvert u_{10}\rvert\delta\right),\qquad && v_{2}(t) = e^{-\lambda (\tau - t)} \left[v_{2\tau} + O\left(\delta^{2}\right)\right].
\end{aligned}
\end{equation}
\end{myrem}

\begin{mylem}[Case $\lambda_{1} <\lambda_{2} < 2\lambda_{1}$]\label{flowlemma}
There exists \(M > 0\) such that for any sufficiently small \(\delta > 0\), and any \(u_{10}\), \(u_{20}\), \(v_{1\tau}\) and \(v_{2\tau}\), where \(\max \lbrace \lvert u_{10}\rvert, \lvert u_{20}\rvert, \lvert v_{1\tau}\rvert, \lvert v_{2\tau}\rvert\rbrace \leq \delta\), the solution \(\left(u\left(t\right), v\left(t\right)\right)\) of system (\ref{eq19000}) that satisfies boundary condition (\ref{eq800020}) can be written as
\begin{equation}\label{eq67900}
\begin{aligned}
u_{1}(t) =& e^{-\lambda_{1} t}u_{10} + \xi_{1}\left(x\right),\qquad
&& u_{2}(t) = e^{-\lambda_{2} t}u_{20} + \xi_{2}\left(x\right),\\
v_{1}(t) =& e^{-\lambda_{1}\left(\tau - t\right)} v_{1\tau} + \zeta_{1}\left(x\right),\qquad
&& v_{2}(t) = e^{-\lambda_{2}\left(\tau-t\right)} v_{2\tau} +\zeta_{2}\left(x\right),
\end{aligned}
\end{equation}
where \(x = \left(t, \tau, u_{10}, u_{20}, v_{1\tau}, v_{2\tau}\right)\), \(t\in\left[0, \tau\right]\), and
\begin{equation*}
\begin{aligned}
\lvert \xi_{1} \rvert \leq & M\left[e^{-\lambda_{1}t}\delta \lvert u_{10}\rvert + e^{-\lambda_{1}\left(\tau - t\right) - \lambda_{2}t} \delta\lvert v_{1\tau}\rvert\right],
&& \lvert \xi_{2} \rvert \leq M e^{-\lambda_{2}t}\delta^{2},\\
\lvert \zeta_{1} \rvert \leq & M \left[e^{-\lambda_{1}\left(\tau - t\right)} \delta \lvert v_{1\tau}\rvert + e^{-\lambda_{2}\left(\tau - t\right) - \lambda_{1}t} \delta\lvert u_{10}\rvert\right],\qquad
&& \lvert \zeta_{2} \rvert \leq M e^{-\lambda_{2}\left(\tau-t\right)}\delta^{2}.
\end{aligned}
\end{equation*}
\end{mylem}

\begin{myrem}
For simplicity, we can write the solution given by Lemma \ref{flowlemma} as
\begin{equation}\label{eqiomljd784987qpqpla874}
\begin{aligned}
u_{1}(t) =& e^{-\lambda_{1}t}u_{10}\left[1 + O\left(\delta\right)\right] + e^{-\lambda_{1}\left(\tau - t\right) - \lambda_{2}t} O\left(\delta \lvert v_{1\tau}\rvert\right),\quad
&& u_{2}(t) = e^{-\lambda_{2}t}\left[u_{20} + O\left(\delta^{2}\right)\right],\\
v_{1}(t) =& e^{-\lambda_{1}\left(\tau - t\right)} v_{1\tau} \left[1 + O\left(\delta\right)\right] + e^{-\lambda_{2}\left(\tau - t\right) - \lambda_{1}t} O\left(\delta \lvert u_{10}\rvert\right),\quad
&& v_{2}(t) = e^{-\lambda_{2}\left(\tau-t\right)}\left[v_{2\tau} + O\left(\delta^{2}\right)\right].
\end{aligned}
\end{equation}
\end{myrem}

\begin{mylem}[Case $\lambda_{2} = 2\lambda_{1}$]\label{flowlemmaResonant2}
Let $\lambda = \lambda_{1} = \frac{\lambda_{2}}{2}$. There exists \(M > 0\) such that for any sufficiently small \(\delta > 0\), and any \(u_{10}\), \(u_{20}\), \(v_{1\tau}\) and \(v_{2\tau}\), where \(\max \lbrace \lvert u_{10}\rvert, \lvert u_{20}\rvert, \lvert v_{1\tau}\rvert, \lvert v_{2\tau}\rvert\rbrace \leq \delta\), the solution \(\left(u\left(t\right), v\left(t\right)\right)\) of system (\ref{eq19000}) that satisfies boundary condition (\ref{eq800020}) can be written in the form (\ref{eq67900}), where \(t\in\left[0, \tau\right]\) and
\begin{equation*}
\begin{aligned}
\lvert \xi_{1} \rvert \leq & M\left[e^{-\lambda t} \lvert u_{10}\rvert\delta + t e^{-\lambda(\tau + t)}\lvert v_{1\tau}\rvert\delta\right],
&& \lvert \xi_{2} \rvert \leq M \left[t e^{-2\lambda t} u_{10}^{2} + e^{-2\lambda t} \delta^{2}\right],\\
\lvert \zeta_{1} \rvert \leq & M \left[e^{-\lambda(\tau - t)} \lvert v_{1\tau}\rvert\delta + \left(\tau - t\right) e^{-\lambda (2\tau - t)}\lvert u_{10}\rvert \delta\right],\quad
&& \lvert \zeta_{2} \rvert \leq M \left[\left(\tau - t\right) e^{-2\lambda(\tau - t)} v_{1\tau}^{2} + e^{-2\lambda(\tau - t)} \delta^{2}\right].
\end{aligned}
\end{equation*}
\end{mylem}

\begin{myrem}
For simplicity, we can write the solution given by Lemma \ref{flowlemmaResonant2} as
\begin{equation}\label{equ13rhuef8yreiosueygs}
\begin{aligned}
u_{1}(t) =& e^{-\lambda t} u_{10}\left[1 + O(\delta)\right] + t e^{-\lambda(\tau + t)} O\left(\lvert v_{1\tau}\rvert\delta\right),\\
u_{2}(t) =& e^{-2\lambda t} u_{20} + t e^{-2\lambda t} O\left(\lvert u_{10}\rvert^{2}\right) + e^{-2\lambda t} O\left(\delta^{2}\right),\\
v_{1}(t) =& e^{-\lambda(\tau - t)} v_{1\tau} \left[1 + O\left(\delta\right)\right] + \left(\tau - t\right) e^{-\lambda (2\tau - t)} O\left(\lvert u_{10}\rvert \delta\right),\\
v_{2}(t) =& e^{-2\lambda(\tau - t)} v_{2\tau} + \left(\tau - t\right) e^{-2\lambda(\tau - t)} O\left(v_{1\tau}^{2}\right) + e^{-2\lambda(\tau - t)} O\left(\delta^{2}\right).
\end{aligned}
\end{equation}
\end{myrem}

\begin{mylem}[Case $\lambda_{2} > 2\lambda_{1}$]\label{flowlemma0}
There exists \(M > 0\) such that for any sufficiently small \(\delta > 0\), and any \(u_{10}\), \(u_{20}\), \(v_{1\tau}\) and \(v_{2\tau}\), where \(\max \lbrace \lvert u_{10}\rvert, \lvert u_{20}\rvert, \lvert v_{1\tau}\rvert, \lvert v_{2\tau}\rvert\rbrace \leq \delta\), the solution \(\left(u\left(t\right), v\left(t\right)\right)\) of system (\ref{eq23000}) that satisfies boundary condition (\ref{eq800020}) can be written in the form (\ref{eq67900}), where \(t\in\left[0, \tau\right]\) and
\begin{equation*}
\begin{aligned}
\lvert \xi_{1} \rvert \leq & M\left[e^{-\lambda_{1}t}\delta \lvert u_{10}\rvert + e^{-\lambda_{1}\left(\tau + t\right)} \delta\lvert v_{1\tau}\rvert\right],
&& \lvert \xi_{2} \rvert \leq M e^{-\lambda_{2}t}\delta^{2},\\
\lvert \zeta_{1} \rvert \leq & M \left[e^{-\lambda_{1}\left(\tau - t\right)} \delta \lvert v_{1\tau}\rvert + e^{-\lambda_{1}\left(2\tau + t\right)} \delta\lvert u_{10}\rvert\right],\qquad
&& \lvert \zeta_{2} \rvert \leq M e^{-\lambda_{2}\left(\tau-t\right)}\delta^{2}.
\end{aligned}
\end{equation*}
\end{mylem}

\begin{myrem}
For simplicity, we can write the solution given by Lemma \ref{flowlemma0} as
\begin{equation}\label{equos83km3odnin8b83xbv}
\begin{aligned}
u_{1}\left(t\right) &= e^{-\lambda_{1}t}u_{10}\left[1 + O\left(\delta\right)\right] + e^{-\lambda_{1}\left(\tau + t\right)} O\left(\delta \lvert v_{1\tau}\rvert\right),\quad
&& u_{2}\left(t\right) = e^{-\lambda_{2}t}\left[u_{20} + O\left(\delta^{2}\right)\right],\\
v_{1}\left(t\right) &= e^{-\lambda_{1}\left(\tau-t\right)}v_{1\tau}\left[1 + O\left(\delta\right)\right] + e^{-\lambda_{1}\left(2\tau - t\right)}O\left(\delta \lvert u_{10}\rvert\right),\quad
&& v_{2}\left(t\right) = e^{-\lambda_{2}\left(\tau-t\right)}\left[v_{2\tau} + O\left(\delta^{2}\right)\right].
\end{aligned}
\end{equation}
\end{myrem}

Lemmas \ref{flowlemma} and \ref{flowlemma0} are proved in \cite{BakraniPhDthesis}, Section 3.3.2. We prove Lemmas \ref{flowlemmaresonant} and \ref{flowlemmaResonant2} below.

\begin{proof}[Proof of Lemma \ref{flowlemmaresonant}]
Let \(x=(u, v)\), \(x\left(t\right) = \left(u\left(t\right), v\left(t\right)\right)\), \(x\left(s\right) = \left(u\left(s\right), v\left(s\right)\right)\), and $\Omega$ be the set of all the points $(u,v)\in\mathbb{R}^{4}$ such that $ \max\{\lvert u_1\rvert, \lvert u_2\rvert, \lvert v_1\rvert, \lvert v_2\rvert\} \leq \delta$. Consider the set
\begin{equation}\label{eqmjhscfnhf6ba7i6w38nrgfjsyg}
\begin{aligned}
A = \big\{x\left(t\right): \quad &\lvert u_{1}(t)\rvert \leq 2 e^{-\lambda t} \lvert u_{10}\rvert + e^{-\lambda\tau} \lvert v_{1\tau}\rvert,\qquad \lvert u_{2}(t)\rvert \leq 2e^{-\lambda t}\delta,\\ &\lvert v_{1}(t)\rvert \leq 2 e^{-\lambda (\tau - t)} \lvert v_{1\tau}\rvert + e^{-\lambda \tau} \lvert u_{10}\rvert,\qquad \lvert v_{2}(t)\rvert \leq 2e^{-\lambda\left(\tau-t\right)}\delta \big\},
\end{aligned}
\end{equation}
where $x: t\mapsto x(t) \in \Omega$ is any continuous function defined for $t\in [0, \tau]$. We show that $A$ is invariant under the integral operator \(\mathfrak{T}\), i.e. \(\mathfrak{T}\left(A\right) \subseteq A\). Then, since $x(t) \equiv (0,0,0,0)$ lies in $A$, and $A$ is closed in $\mathcal{C}^{0}([0,\tau], \Omega)$, Shilnikov's theorem on boundary value problems implies that the fixed point of the integral operator $\mathfrak{T}$ which is indeed the solution \(\left(u\left(t\right), v\left(t\right)\right)\) of system (\ref{eq74wt366cuw6gt5uw6v01092}) lies in $A$ too.

To show the invariance of $A$, recast system (\ref{eq74wt366cuw6gt5uw6v01092}) into the form (\ref{eqbwiebro4ibyw43yyyiw0290}), where
\begin{equation}\label{eqgyyty736q487eg6347ctsg}
\begin{aligned}
F_{1}\left(x\right) =& \mathtt{f}_{11}\left(x\right) u_{1}^{2} + \mathtt{f}_{12}\left(x\right) u_{1} u_{2} + \mathtt{f}_{13}\left(x\right) u_{1} v_{1} + \mathtt{f}_{14}\left(x\right) u_{1} v_{2} + \mathtt{f}_{15}\left(x\right) u_{2} v_{1},\\
F_{2}\left(x\right) =& \mathtt{f}_{21}\left(x\right) u_{1}^{2} + \mathtt{f}_{22}\left(x\right) u_{1} u_{2} + \mathtt{f}_{23}\left(x\right) u_{1} v_{1} + \mathtt{f}_{24}\left(x\right) u_{1} v_{2} + \mathtt{f}_{25}\left(x\right) u_{2}^{2}\\
&+ \mathtt{f}_{26}\left(x\right) u_{2} v_{1} + \mathtt{f}_{27}\left(x\right) u_{2} v_{2},\\
G_{1}\left(x\right) =& \mathtt{g}_{11}\left(x\right) v_{1}^{2} + \mathtt{g}_{12}\left(x\right) v_{1} v_{2} + \mathtt{g}_{13}\left(x\right) v_{1} u_{1} + \mathtt{g}_{14}\left(x\right) v_{1} u_{2} + \mathtt{g}_{15}\left(x\right) v_{2} u_{1},\\
G_{2}\left(x\right) =& \mathtt{g}_{21}\left(x\right) v_{1}^{2} + \mathtt{g}_{22}\left(x\right) v_{1} v_{2} + \mathtt{g}_{23}\left(x\right) v_{1} u_{1} + \mathtt{g}_{24}\left(x\right) v_{1} u_{2} + \mathtt{g}_{25}\left(x\right) v_{2}^{2}\\
&+ \mathtt{g}_{26}\left(x\right) v_{2} u_{1} + \mathtt{g}_{27}\left(x\right) v_{2} u_{2},
\end{aligned}
\end{equation}
for some smooth functions \(\mathtt{f}_{ij}\) and \(\mathtt{g}_{ij}\) defined on $\Omega$. Let $M_1 := \max\{\lvert \mathtt{f}_{ij}(u,v)\rvert, \lvert \mathtt{g}_{ij}(u,v)\rvert\}$. Since $\Omega$ is compact, $M_1$ is well-defined. By (\ref{eqmjhscfnhf6ba7i6w38nrgfjsyg}) and (\ref{eqgyyty736q487eg6347ctsg}), for any \(x(t)\in A\), we have
\begin{equation*}
\begin{aligned}
\left\vert F_{1}\left(x\left(t\right)\right)\right\vert & \leq M_1 M_2 \left(e^{- 2\lambda t}\lvert u_{10}\rvert\delta + e^{-\lambda\tau} \lvert u_{10}\rvert\delta + e^{-\lambda\tau} \lvert v_{1\tau}\rvert\delta\right),\\
%%%%%%
\left\vert F_{2}\left(x\left(t\right)\right)\right\vert & \leq M_1 M_2 \left(e^{- \lambda \tau} \delta^{2} + e^{-2\lambda t} \delta^{2}\right),\\
%%%%%%
\left\vert G_{1}\left(x\left(t\right)\right)\right\vert & \leq M_1 M_2 \left(e^{- 2\lambda (\tau -t)}\lvert v_{1\tau}\rvert\delta + e^{-\lambda\tau} \lvert v_{1\tau}\rvert\delta + e^{-\lambda\tau} \lvert u_{10}\rvert\delta\right),\\
\left\vert G_{2}\left(x\left(t\right)\right)\right\vert & \leq M_1 M_2 \left(e^{-\lambda \tau} \delta^{2} + e^{-2\lambda t}\delta^{2}\right)
\end{aligned}
\end{equation*}
for some $M_2 > 0$ (see \cite{BakraniPhDthesis} for the calculation of such a constant). Define $M := 2 M_1 M_2 \lambda^{-1}$. Then,
\begin{equation*}
\begin{aligned}
\left\vert \overline{u}_{1}\left(t\right) - e^{-\lambda t} u_{10}\right\vert \leq & \int_{0}^{t} e^{\lambda(s-t)} \left\vert F_{1}\left(x\left(s\right)\right)\right\vert ds
\leq M_1 M_2 \int_{0}^{t} e^{\lambda(s-t)} \Big(e^{- 2\lambda s}\lvert u_{10}\rvert + e^{-\lambda\tau} \lvert u_{10}\rvert\\
&+ e^{-\lambda\tau} \lvert v_{1\tau}\rvert\Big)\delta ds \leq M \left(e^{-\lambda t}\lvert u_{10}\rvert + e^{-\lambda\tau} \lvert v_{1\tau}\rvert\right) \delta,\\
%%%%%
%%%%%
\left\vert \overline{u}_{2}\left(t\right) - e^{-\lambda t} u_{20}\right\vert \leq & \int_{0}^{t} e^{\lambda(s-t)} \left\vert F_{2}\left(x\left(s\right)\right)\right\vert ds
\leq M_1 M_2 \int_{0}^{t} e^{\lambda(s-t)}\left(e^{-\lambda \tau} + e^{-2\lambda s}\right) \delta^2 ds
\leq M e^{-\lambda t} \delta^2,
\end{aligned}
\end{equation*}
\begin{equation*}
\begin{aligned}
\left\vert \overline{v}_{1}\left(t\right) - e^{-\lambda\left(\tau-t\right)} v_{1\tau}\right\vert
\leq & \int_{t}^{\tau} e^{\lambda(t-s)} \lvert G_{1}(x(s))\rvert ds
\leq M_1 M_2 \int_{t}^{\tau} e^{\lambda(t-s)}\Big(e^{- 2\lambda (\tau - s)}\lvert v_{1\tau}\rvert + e^{-\lambda\tau} \lvert v_{1\tau}\rvert\\
& + e^{-\lambda\tau} \lvert u_{10}\rvert\Big) \delta ds
\leq M \left(e^{-\lambda (\tau -t)}\lvert v_{1\tau}\rvert + e^{-\lambda\tau} \lvert u_{10}\rvert\right) \delta,\\
%%%%%%
%%%%%%
\left\vert \overline{v}_{2}\left(t\right) - e^{-\lambda (\tau - t)} v_{2\tau}\right\vert \leq & \int_{t}^{\tau} e^{\lambda(t-s)} \left\vert G_{2}(x(s))\right\vert ds
\leq M_1 M_2 \int_{t}^{\tau} e^{\lambda(t-s)}\left(e^{-\lambda \tau} + e^{-2\lambda (\tau - s)}\right) \delta^2 ds\\
\leq & M e^{-\lambda (\tau - t)} \delta^2.
\end{aligned}
\end{equation*}
The invariance of $A$ is derived from these relations if $\delta$ is chosen sufficiently small, i.e. $\delta < \frac{1}{M}$. Meanwhile, we have shown that the image of any element of $A$ under the integral operator $\mathfrak{T}$ is of the form (\ref{eq879ygo7f6gvjhferg}). However, since the solution \(\left(u\left(t\right), v\left(t\right)\right)\) of system (\ref{eq74wt366cuw6gt5uw6v01092}) is the fixed point of $\mathfrak{T}$, we have that it must be of the form (\ref{eq879ygo7f6gvjhferg}) too. This proves Lemma \ref{flowlemmaresonant}.
\end{proof}

\begin{proof}[Proof of Lemma \ref{flowlemmaResonant2}]
We follow the same notations and procedure as in the proof of Lemma \ref{flowlemmaresonant}. Let \(\delta > 0\) be small, and consider the set
\begin{equation}\label{eq8b9q20b3209bql02rv48t6v476}
\begin{aligned}
A = \Big\{x(t):\quad & \lvert u_{1}(t)\rvert \leq 2e^{-\lambda t}\lvert u_{10}\rvert + e^{-\lambda\left(\tau + t\right)}\lvert v_{1\tau}\rvert,\qquad \lvert u_{2}(t)\rvert \leq 2 e^{-2\lambda t}\delta,\\
& \lvert v_{1}(t)\rvert \leq 2e^{-\lambda \left(\tau-t\right)}\lvert v_{1\tau}\rvert + e^{-\lambda\left(2\tau - t\right)} \lvert u_{10}\rvert,\qquad \lvert v_{2}(t)\rvert \leq 2 e^{-2\lambda\left(\tau-t\right)}\delta\Big\},
\end{aligned}
\end{equation}
where \(\left(u_{1}(t), u_{2}(t), v_{1}(t), v_{2}(t)\right)\) is any continuous function defined on \(\Omega\) for \(t \in \left[0, \tau\right]\). System (\ref{eq19000}) satisfies (\ref{eq20000}) and (\ref{eq24500}), and so we can write it in the form (\ref{eqbwiebro4ibyw43yyyiw0290}), where
\begin{equation}\label{eqksey48b7itybi47vw0901wn}
\begin{aligned}
F_{1}\left(u,v\right) &= \mathtt{f}_{11}\left(u,v\right) u_{1}^{2} + \mathtt{f}_{12}\left(u,v\right) u_{1} u_{2} + \mathtt{f}_{13}\left(u,v\right) v_{1} u_{2},\\
F_{2}\left(u,v\right) &= \mathtt{f}_{21}\left(u,v\right) u_{1}^{2} + \mathtt{f}_{22}\left(u,v\right) u_{1} u_{2} + \mathtt{f}_{23}\left(u,v\right) u_{2}^{2},\\
G_{1}\left(u,v\right) &= \mathtt{g}_{11}\left(u,v\right) v_{1}^{2} + \mathtt{g}_{12}\left(u,v\right) v_{1} v_{2} + \mathtt{g}_{13}\left(u,v\right) u_{1} v_{2},\\
G_{2}\left(u,v\right) &= \mathtt{g}_{21}\left(u,v\right) v_{1}^{2} + \mathtt{g}_{22}\left(u,v\right) v_{1} v_{2} + \mathtt{g}_{23}\left(u,v\right) v_{2}^{2},
\end{aligned}
\end{equation}
for some continuous functions \(\mathtt{f}_{ij}\) and \(\mathtt{g}_{ij}\). Thus, for any $x(t)\in A$, we have
\begin{equation*}
\begin{aligned}
\lvert F_{1}\left(x\left(t\right)\right)\rvert &\leq M_{1} M_{2}\left(e^{-2\lambda t} u_{10}^{2} + e^{-3\lambda t} \lvert u_{10}\rvert\delta + e^{-\lambda (\tau + t)} \lvert v_{1\tau}\rvert\delta\right),\\
%%%%%%%%%%%
\lvert F_{2}\left(x\left(t\right)\right)\rvert &\leq M_{1} M_{2}\left(e^{-2\lambda t} u_{10}^{2} + e^{-3\lambda t} \lvert u_{10}\rvert\delta + e^{-4\lambda t}\delta^{2}\right),\\
%%%%%%%%%%%
\lvert G_{1}\left(x\left(t\right)\right)\rvert &\leq M_{1} M_{2}\left(e^{-2\lambda(\tau - t)} v_{1\tau}^{2} + e^{-3\lambda(\tau - t)} \lvert v_{1\tau}\rvert\delta + e^{-\lambda(2\tau - t)} \lvert u_{10}\rvert\delta\right),\\
%%%%%%%%%%%
\lvert G_{2}\left(x\left(t\right)\right)\rvert &\leq M_{1} M_{2}\left(e^{-2\lambda(\tau - t)} v_{1\tau}^{2} + e^{-3\lambda(\tau - t)} \lvert v_{1\tau}\rvert\delta + e^{-4\lambda(\tau - t)} \delta^{2}\right)\end{aligned}
\end{equation*}
for some positive $M_{1}$ and $M_{2}$. Then, there exists an $M>0$ such that
\begin{equation*}
\begin{aligned}
\left\vert \overline{u}_{1}\left(t\right) - e^{-\lambda t} u_{10}\right\vert \leq & \int_{0}^{t} e^{\lambda(s-t)} \left\vert F_{1}\left(x\left(s\right)\right)\right\vert ds
\leq M_1 M_2 \int_{0}^{t} e^{\lambda(s-t)} \Big(e^{-2\lambda s} u_{10}^{2} + e^{-3\lambda s} \lvert u_{10}\rvert\delta\\
&+ e^{-\lambda (\tau + s)} \lvert v_{1\tau}\rvert\delta\Big) ds
\leq M\left(e^{-\lambda t} \lvert u_{10}\rvert\delta + t e^{-\lambda(\tau + t)}\lvert v_{1\tau}\rvert\delta\right),\\
%%%%%%%%%%
\left\vert \overline{u}_{2}\left(t\right) - e^{-2\lambda t} u_{20}\right\vert \leq & \int_{0}^{t} e^{2\lambda(s-t)} \left\vert F_{2}\left(x\left(s\right)\right)\right\vert ds
\leq M_1 M_2 \int_{0}^{t} e^{2\lambda(s-t)} \Big(e^{-2\lambda s} u_{10}^{2} + e^{-3\lambda s} \lvert u_{10}\rvert\delta\\
&+ e^{-4\lambda s}\delta^{2}\Big) ds
\leq M \left(t e^{-2\lambda t} u_{10}^{2} + e^{-2\lambda t} \delta^{2}\right),
\end{aligned}
\end{equation*}
\begin{equation*}
\begin{aligned}
\left\vert \overline{v}_{1}\left(t\right) - e^{-\lambda\left(\tau-t\right)} v_{1\tau}\right\vert
\leq & \int_{t}^{\tau} e^{\lambda(t-s)} \lvert G_{1}\left(x\left(s\right)\right)\rvert ds
\leq M_1 M_2 \int_{t}^{\tau} e^{\lambda(t-s)}\Big(e^{-2\lambda(\tau - s)} v_{1\tau}^{2} + e^{-3\lambda(\tau - s)} \lvert v_{1\tau}\rvert\delta\\
&+ e^{-\lambda(2\tau - s)} \lvert u_{10}\rvert\delta\Big) ds
\leq M \left(e^{-\lambda(\tau - t)} \lvert v_{1\tau}\rvert\delta + \left(\tau - t\right) e^{-\lambda (2\tau - t)}\lvert u_{10}\rvert \delta\right),\\
%%%%%%
%%%%%%
\left\vert \overline{v}_{2}\left(t\right) - e^{-2\lambda (\tau - t)} v_{2\tau}\right\vert \leq & \int_{t}^{\tau} e^{2\lambda(t-s)} \left\vert G_{2}\left(x\left(s\right)\right)\right\vert ds
\leq M_1 M_2 \int_{t}^{\tau} e^{2\lambda(t-s)}\Big(e^{-2\lambda(\tau - s)} v_{1\tau}^{2} + e^{-3\lambda(\tau - s)} \lvert v_{1\tau}\rvert\delta\\
&+ e^{-4\lambda(\tau - s)} \delta^{2}\Big) ds \leq M \left(\left(\tau - t\right) e^{-2\lambda(\tau - t)} v_{1\tau}^{2} + e^{-2\lambda(\tau - t)} \delta^{2}\right).
\end{aligned}
\end{equation*}
The rest of the proof is the same as the proof of Lemma \ref{flowlemmaresonant}. This ends the proof of Lemma \ref{flowlemmaResonant2}.
\end{proof}

\section{Proof of Theorem \ref{Thm_Single_Loop}}\label{Sec_Proof_Thm_Single_Loop}

We prove Theorem \ref{Thm_Single_Loop} in this section. The proof is based on the study of the Poincar\'e map along $\Gamma$. Recall $\Sigma^{\mathrm{in}}$, $\Sigma^{\mathrm{out}}$, $\Pi^{\mathrm{in}}(h)$ and $\Pi^{\mathrm{out}}(h)$ from the Introduction (see Figure \ref{Fig_three_dim_cross_sections}). Fix an $h$ such that $\lvert h\rvert \leq h_0$ for some sufficiently small $h_{0} > 0$. Consider an orbit near $\Gamma$ that starts from $\Pi^{\mathrm{in}}(h)$ at a point $(u_{10}, \delta, v_{10}, v_{20})$, goes along $\Gamma$ and after passing time $\tau >0$ intersects $\Pi^{\mathrm{out}}(h)$ at a point $(u_{1\tau}, u_{2\tau}, v_{1\tau}, \delta)$. This orbit keeps going along $\Gamma$ until it returns to $\Pi^{\mathrm{in}}(h)$ at a point $(\overline{u}_{10}, \delta, \overline{v}_{10}, \overline{v}_{20})$ (not every orbit starting from $\Pi^{\mathrm{in}}(h)$ close to $\Gamma$ takes such a journey; some orbits may go along the negative side of the unstable manifold of $O$). Such orbits define a Poincar\'e map (first return map) $T_h$ along $\Gamma$. As this map appears frequently in this paper, we drop the subscript $h$ from $T_h$ to avoid unnecessary complications. To study this map, we first endow the cross-sections with lower dimensional coordinates.
\begin{mylem}\label{Coordinates_Section_Lemma}
We can choose $(u_{1}, v_{1})$-coordinates on each of the cross-sections $\Pi^{\mathrm{in}}(h)$ and $\Pi^{\mathrm{out}}(h)$. More precisely, $(u_{2}, v_{2})$ is uniquely determined by $(u_{1}, v_{1})$ for any $(u_{1}, u_{2}, v_{1}, v_{2})$ on $\Pi^{\mathrm{in}}(h)$ or $\Pi^{\mathrm{out}}(h)$.
\end{mylem}
\begin{proof}
Let $H$ be the first integral of any of systems (\ref{eq74wt366cuw6gt5uw6v01092}), (\ref{eq19000}) and (\ref{eq23000}). Then, $\frac{\partial H}{\partial u_{2}}(0,0,0,\delta)$, $\frac{\partial H}{\partial v_{2}}(0,\delta,0,0)$, $\frac{\partial H}{\partial u_{2}}(0,0,0,-\delta)$ and $\frac{\partial H}{\partial v_{2}}(0, -\delta,0,0)$ are nonzero (see \cite{BakraniPhDthesis}, Section 3.2.1). Thus, having the relation $H=h$, the implicit function theorem states that $(u_{2}, v_{2})$ is a smooth function of $(u_{1}, v_{1})$ near $(u_{1}, v_{1}) = (0,0)$, as desired.
\end{proof}

\begin{mynotation}
We denote the point $(0,0)$ on $\Pi^{\mathrm{in}}(h)$ (resp. $\Pi^{\mathrm{out}}(h)$) by $M^{\mathrm{in}}(h)$ (resp. $M^{\mathrm{out}}(h)$).
\end{mynotation}

The periodic orbit $L_h$ intersects both $\Pi^{\mathrm{in}}(h)$ and $\Pi^{\mathrm{out}}(h)$ at $(u_{1}, v_{1}) = (0,0)$ as it entirely lies in the plane $\{u_{1} = v_{1} = 0\}$ (it intersect the sections at $M^{\mathrm{in}}(h)$ and $M^{\mathrm{out}}(h)$) (see Figure \ref{Fig_three_dim_cross_sections}).

Considering $(u_{1}, v_{1})$-coordinates on the cross-sections, we define the Poincar\'e map $T: \mathcal{D}\subset \Pi^{\mathrm{in}}(h) \rightarrow \Pi^{\mathrm{in}}(h)$, where $\mathcal{D}$ is the domain of the map, by
\begin{equation*}
(u_{10}, v_{10}) \mapsto (\overline{u}_{10}, \overline{v}_{10}).
\end{equation*}
An orbit starting from a point on $\Pi^{\mathrm{out}}(h)$ close to $M^{\mathrm{out}}(h)$ goes along $L_{h}$ until it intersects $\Pi^{\mathrm{in}}(h)$ at a point close to $M^{\mathrm{in}}(h)$. Such orbits define a diffeomorphism $T^{\mathrm{glo}}$, call it global map, from a neighborhood of $M^{\mathrm{out}}(h)$ in $\Pi^{\mathrm{out}}(h)$ to $\Pi^{\mathrm{in}}(h)$. In particular, we have that $T^{\mathrm{glo}}(u_{1\tau}, v_{1\tau}) = (\overline{u}_{10}, \overline{v}_{10})$. Then, by defining a local map $T^{\mathrm{loc}}: (u_{10}, v_{10})\mapsto (u_{1\tau}, v_{1\tau})\in \Pi^{\mathrm{out}}(h)$ for $(u_{10}, v_{10})\in\mathcal{D}$, we can write
\begin{equation}\label{eq_composition_of_global_and_local_maps}
T = T^{\mathrm{glo}} \circ T^{\mathrm{loc}}.
\end{equation}

The global map $T^{\mathrm{glo}}$ is a diffeomorphism and well-approximated by its Taylor expansion at $M^{\mathrm{out}}(h) \in \Pi^{\mathrm{out}}(h)$, i.e.
\begin{equation}\label{eq_Taylor_global_map}
T^{\text{glo}}\left(\begin{matrix}
u_1 \\ v_1
\end{matrix}\right) = \left(\begin{matrix}
a(h) & b(h)\\ c(h) & d(h)
\end{matrix}\right) \left(\begin{matrix}
u_1 \\ v_1
\end{matrix}\right) + \left(\begin{matrix}
o\left(\lvert u_{1}\rvert, \lvert v_{1}\rvert\right)\\
o\left(\lvert u_{1}\rvert, \lvert v_{1}\rvert\right)
\end{matrix}\right).
\end{equation}
Here, $a(h)$, $b(h)$, $c(h)$ and $d(h)$ are smooth real-valued functions of $h$. Evaluating these functions at $h=0$ gives the constants $a$, $b$, $c$ and $d$ in (\ref{eq63000}), which are indeed the constants used in \cite{Bakrani2022JDE} for the study of the dynamics in the zero level set.

The domain $\mathcal{D}$ is defined by the orbits that start from $\Pi^{\mathrm{in}}(h)$, go along $\Gamma$ (or equivalently, go along $L_h$ for $h$ sufficiently close to zero), intersect $\Pi^{\mathrm{out}}(h)$, and then $\Pi^{\mathrm{in}}(h)$ again. The following lemma follows from the smoothness of the system and the structure of the reduced system to the invariant plane (see also Figures \ref{Periodic_orbit_single_homoclinic_case} and \ref{Periodic_orbits_double_homoclinic_case}).

\begin{mylem}\label{Lem76tgitiewr}
If $(0,0)\in\mathcal{D}$, then $\mathcal{D}$ contains a small open neigborhood of $(0,0)$ in $\Pi^{\mathrm{in}}(h)$. If $(0,0)\notin\mathcal{D}$, no point in some small open neigborhood of $(0,0)$ in $\Pi^{\mathrm{in}}(h)$ lies in $\mathcal{D}$.
\end{mylem}

Note also that in the case $(0,0)\in\mathcal{D}$, the point $(0,0)$ is nothing but the intersection of the periodic orbit $L_{h}$ and $\Pi^{\mathrm{in}}_{h}$, i.e. $M^{\mathrm{in}}(h)$ is a fixed point of $T$.

To study the domain $\mathcal{D}$, we only focus on the orbits that start from the $\epsilon$-neighborhood of $(0,0)$ in $\Pi^{\mathrm{in}}(h)$ and intersect the $\epsilon$-neighborhood of $(0,0)$ in $\Pi^{\mathrm{out}}(h)$ for some sufficiently small $0 < \epsilon < \delta$, i.e.
\begin{equation}\label{eqDomain}
\mathcal{D} = \{(u_{10}, v_{10})\in \Pi^{\mathrm{in}}(h): \quad \|(u_{10}, v_{10})\|\leq \epsilon,\quad (u_{1\tau}, v_{1\tau})\in\Pi^{\mathrm{out}}(h) \quad \mathrm{and}\quad \|(u_{1\tau}, v_{1\tau})\|\leq \epsilon\}.
\end{equation}
Let \(B_{\epsilon}\) be the open \(\epsilon\)-ball in \(\Pi^{\mathrm{in}}(h)\) centered at \((0,0)\). Fix a sufficiently large \(m>1\), and define
\begin{equation}\label{eqyujnj68999jhhhgg34}
\begin{gathered}
Y_{1} := \left\{\left(u_{10}, v_{10}\right)\in B_{\epsilon}: \lvert v_{10}\rvert < m^{-1}\lvert u_{10}\rvert\right\}\quad\mathrm{and}\quad
Y_{2} := \left\{\left(u_{10}, v_{10}\right)\in B_{\epsilon}: m^{-1}\lvert u_{10}\rvert \leq \lvert v_{10}\rvert\right\}.
\end{gathered}
\end{equation}
Obviously, \(B_{\epsilon} = Y_{1} \cup Y_{2}\). We also define
\begin{equation}\label{eqiuvyseyt74vqi8l3viul181}
\mathcal{D}_{1}:= \mathcal{D}\cap Y_{1}\qquad \mathrm{and}\qquad \mathcal{D}_{2}:= \mathcal{D}\cap Y_{2}.
\end{equation}
Obviously, $\mathcal{D} = \mathcal{D}_{1} \cup \mathcal{D}_{2}$.

The next lemma discusses the dynamics of the Poincar\'e map $T$ on the region $\mathcal{D}_{2}$, and $T^{-1}$ on the region $\mathcal{D}_{1}$. It states that the forward (resp. backward) orbit of a point in $\mathcal{D}_{2}$ (resp. $\mathcal{D}_{1}$) never enters $\mathcal{D}_{1}$ (resp. $\mathcal{D}_{2}$). Moreover, $T$ (resp. $T^{-1}$) expands $\mathcal{D}_{2}$ (resp. $\mathcal{D}_{1}$). Thus, starting in $\mathcal{D}_{2}$ (resp. $\mathcal{D}_{1}$) and iterating $T$ (resp. $T^{-1}$), we stay in $\mathcal{D}_{2}$ (resp. $\mathcal{D}_{1}$) for a finite number of iterations until we leave $\mathcal{D}$.

\begin{mylem}\label{Lem11u67tt6}
Reduce system (\ref{eq300}) near the origin to one of the normal forms (\ref{eq74wt366cuw6gt5uw6v01092}), (\ref{eq19000}) or (\ref{eq23000}).
Let $(\overline{u}_{10}, \overline{v}_{10}) = T(u_{10}, v_{10})$, for some $(u_{10}, v_{10})\in \mathcal{D}$. Then, the following hold.
\begin{enumerate}[(i)]
\item Suppose $(u_{10}, v_{10})\in \mathcal{D}_2$ such that $(\overline{u}_{10}, \overline{v}_{10}) \in B_{\epsilon}$. Then $(\overline{u}_{10}, \overline{v}_{10}) \in Y_{2}$. Furthermore, $(u_{10}, v_{10}) = o\left(\|(\overline{u}_{10}, \overline{v}_{10})\|\right)$ as $(\epsilon, h_{0})\rightarrow (0,0)$. In addition, $\frac{\overline{v}_{10}}{\overline{u}_{10}} \rightarrow \frac{d(h)}{b(h)}$ as $(u_{10}, v_{10})\rightarrow (0,0)$.
\item Suppose $(\overline{u}_{10}, \overline{v}_{10}) \in Y_{1}$ for some $(u_{10}, v_{10}) \in \mathcal{D}$. Further assume $(u_{10}, v_{10}) \in T(\mathcal{D})$ (this assumption is natural as we are interested in the points whose backward orbits under $T$ remain in $\mathcal{D}$). Then, $(u_{10}, v_{10}) \in \mathcal{D}_{1}$. Furthermore, $(\overline{u}_{10}, \overline{v}_{10}) = o\left(\|(u_{10}, v_{10})\|\right)$ as $(\epsilon, h_{0})\rightarrow (0,0)$. In addition, $\frac{v_{10}}{u_{10}} \rightarrow 0$ as $(\overline{u}_{10}, \overline{v}_{10})\rightarrow (0,0)$.
\end{enumerate}
\end{mylem}

The following corollary is a direct consequence of this Lemma.
\begin{mycor}\label{Cor24g76ftyfduwq}
Let $(u_{10}, v_{10})\in\mathcal{D}$. If the forward (resp. backward) iterations of $(u_{10}, v_{10})$ under $T$ lies entirely in $\mathcal{D}$, then $(u_{10}, v_{10})$ must lie in $\mathcal{D}_{1}$ (resp. $\mathcal{D}_{2}$), and we have that $T^{n}(u_{10}, v_{10})$ (resp. $T^{-n}(u_{10}, v_{10})$) converges to $(0,0)$ as $n\rightarrow\infty$.
\end{mycor}

To prove Lemma \ref{Lem11u67tt6}, we first state the following auxiliary lemma. The proof of this lemma is quite technical, and we postpone it to Section \ref{Technical_Section}.

\begin{mylem}\label{Lem568g67eywfg76}
Reduce the system to the corresponding normal form (\ref{eq74wt366cuw6gt5uw6v01092}), (\ref{eq19000}) or (\ref{eq23000}). Let $(u_{10}, v_{10}) \in \mathcal{D}_{2}$, and consider the corresponding $(u_{1\tau}, v_{1\tau}) \in \Pi^{\mathrm{out}}(h)$. Then, as $(h_{0}, \epsilon) \rightarrow (0,0)$, we have $v_{1\tau} = e^{\lambda_{1}\tau} v_{10}\left[1 + O\left(\delta\right)\right]$ and $u_{1\tau} = o\left(\lvert v_{1\tau}\rvert\right)$.
\end{mylem}

\begin{proof}[Proof of Lemma \ref{Lem11u67tt6}]
Consider $(u_{10}, v_{10}) \in \mathcal{D}_{2}$ such that $(\overline{u}_{10}, \overline{v}_{10}) = T(u_{10}, v_{10})\in B_{\epsilon}$. By (\ref{eq_Taylor_global_map}), we obtain
\begin{equation*}
\overline{u}_{10} = \left[a(h) + o\left(1\right)\right] u_{1\tau} + \left[b(h) + o\left(1\right)\right] v_{1\tau},\qquad
\overline{v}_{10} = \left[c(h) + o\left(1\right)\right] u_{1\tau} + \left[d(h) + o\left(1\right)\right] v_{1\tau}.
\end{equation*}
By Lemma \ref{Lem568g67eywfg76}, we have $u_{1\tau} = o\left(\lvert v_{1\tau}\rvert\right)$ as $(h,\epsilon) \rightarrow (0,0)$, and so
\begin{equation*}
\overline{u}_{10} = \left[b(h) + o\left(1\right)\right] v_{1\tau},\qquad
\overline{v}_{10} = \left[d(h) + o\left(1\right)\right] v_{1\tau},
\end{equation*}
which gives $\frac{\overline{v}_{10}}{\overline{u}_{10}} = \frac{d(h)}{b(h)} + o\left(1\right)$. Therefore, if $(u_{10}, v_{10})\in \mathcal{D}_2$ and $T(u_{10}, v_{10})\in B_{\epsilon}$, then $T(u_{10}, v_{10})$ is close to the straight line with the slope $\frac{d(h)}{b(h)}$. Thus, if $m>1$ in (\ref{eqyujnj68999jhhhgg34}) is chosen large enough, we have that $T(u_{10}, v_{10})\in Y_2$. We further have $\frac{\overline{v}_{10}}{\overline{u}_{10}} \rightarrow \frac{d(h)}{b(h)}$ as $(u_{10}, v_{10})\rightarrow (0,0)$.

Having $(u_{10}, v_{10}) \in \mathcal{D}_{2}$ implies $\lvert u_{10}\rvert \leq m \lvert v_{10}\rvert$. Thus, for $(u_{10}, v_{10})\neq (0,0)$, Lemma \ref{Lem568g67eywfg76} yields
\begin{equation*}
\frac{\|(u_{10}, v_{10})\|}{\|(\overline{u}_{10}, \overline{v}_{10})\|} \leq \frac{\sqrt{1 + m^{2}} \lvert v_{10}\rvert}{\lvert v_{1\tau}\rvert \sqrt{[d(h)]^{2} + [b(h)]^{2} + o(1)}} = \frac{\sqrt{1 + m^{2}}}{\left[1 + O\left(\delta\right)\right] \sqrt{[d(h)]^{2} + [b(h)]^{2} + o(1)}} \cdot e^{-\lambda_{1} \tau}.
\end{equation*}
As $(\epsilon, h_{0})\rightarrow (0,0)$, we have $\tau = \tau(u_{10}, v_{10}) \rightarrow \infty$. This implies $\|(u_{10}, v_{10})\| = o\left(\|(\overline{u}_{10}, \overline{v}_{10})\|\right)$, as desired. This finishes the proof of the first part of Lemma \ref{Lem11u67tt6}.

We obtain the statement of the second part of Lemma \ref{Lem11u67tt6} from the first part by a reversion of time. Consider the normal form, and let us call it system I. Take a point $x_{0} \in \Pi^{\mathrm{out}}(h)$, and suppose that its forward orbit intersects $\Pi^{\mathrm{in}}(h)$ at a point $x_{1}$, then intersects $\Pi^{\mathrm{out}}(h)$ at a point $x_{2}$, and then intersects $\Pi^{\mathrm{in}}(h)$ at a point $x_{3}$ (the point $x_{3}$ plays the role of $(\overline{u}_{10}, \overline{v}_{10})$ throughout the proof). In summary (see also Figure \ref{Reversed_time_systems}),
\begin{equation*}
x_{0} \in \Pi^{\mathrm{out}}(h) \longrightarrow x_{1}\in \Pi^{\mathrm{in}}(h) \longrightarrow x_{2}\in \Pi^{\mathrm{out}}(h) \longrightarrow x_{3}\in \Pi^{\mathrm{in}}(h).
\end{equation*}
Let system II be the one that is derived from system I by a reversion of time. For system II, we have
\begin{equation*}
x_{3}\in \Pi^{\mathrm{in}}(h) \longrightarrow x_{2}\in \Pi^{\mathrm{out}}(h) \longrightarrow x_{1}\in \Pi^{\mathrm{in}}(h) \longrightarrow x_{0} \in \Pi^{\mathrm{out}}(h).
\end{equation*}
Let system III be the one which is derived from system II by applying the linear change of coordinates $\left(\widetilde{u}_{1}, \widetilde{u}_{2}, \widetilde{v}_{1}, \widetilde{v}_{2}\right) = \left(v_{1}, v_{2}, u_{1}, u_{2}\right)$. Indeed, system III is derived from system II by exchanging the stable and unstable components. Let $\Pi^{\mathrm{in}}_{\mathrm{III}}(h) := \{\widetilde{u}_{2} = \delta\}\cap\{H=h\}$, $\Pi^{\mathrm{out}}_{\mathrm{III}}(h) := \{\widetilde{v}_{2} = \delta\}\cap\{H=h\}$, and $J := \left(\begin{smallmatrix}
0 & 1\\
1 & 0
\end{smallmatrix}\right)$. Then, for system III, we have
\begin{equation*}
J x_{3}\in \Pi^{\mathrm{out}}_{\mathrm{III}}(h) \longrightarrow
J x_{2}\in \Pi^{\mathrm{in}}_{\mathrm{III}}(h)\longrightarrow
J x_{1}\in \Pi^{\mathrm{out}}_{\mathrm{III}}(h) \longrightarrow
J x_{0}\in \Pi^{\mathrm{in}}_{\mathrm{III}}(h).
\end{equation*}

\begin{figure}[h]
\centering
\begin{subfigure}{0.32\textwidth}
\centering
\includegraphics[scale=0.11]{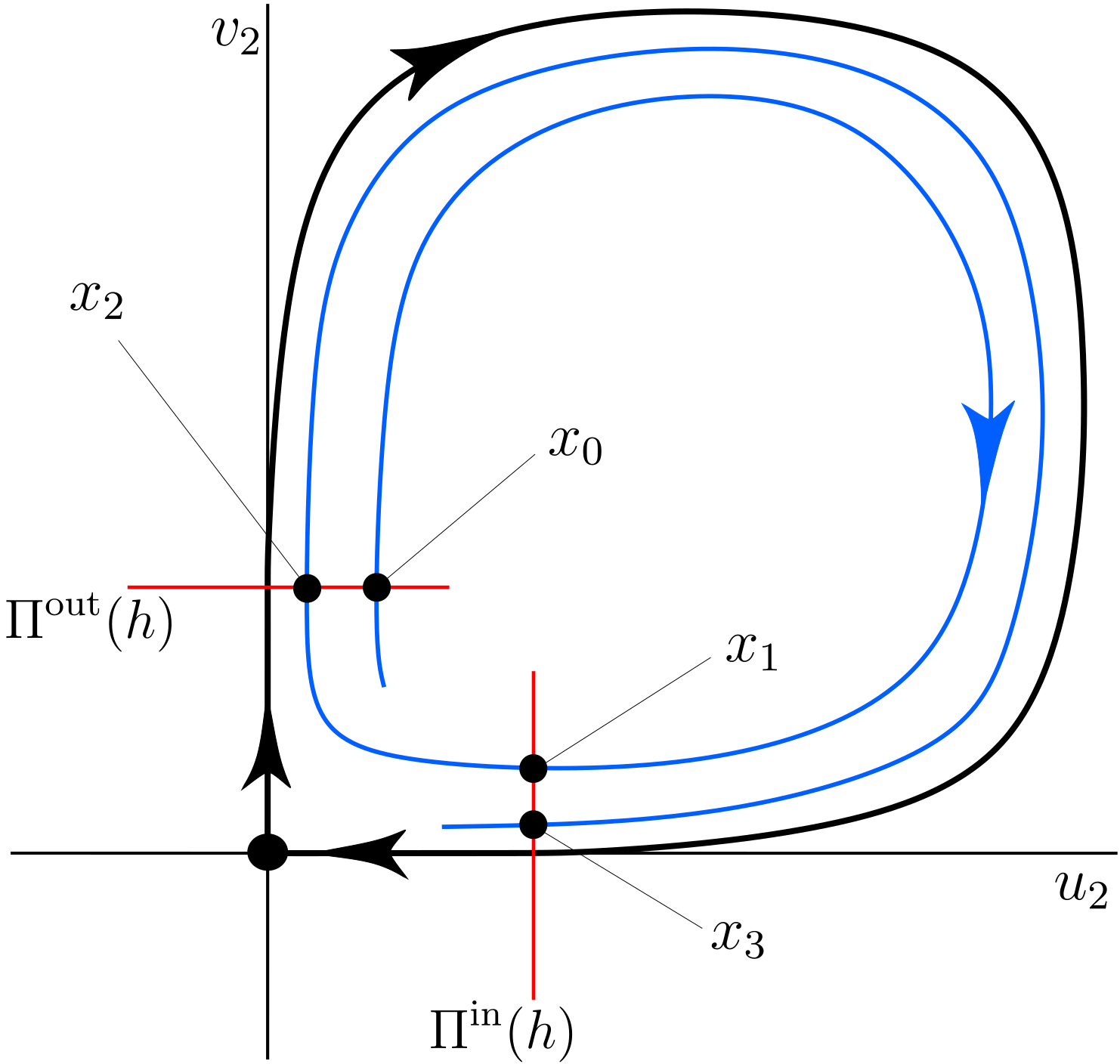}
\caption{System I}
\end{subfigure}
\begin{subfigure}{0.32\textwidth}
\centering
\includegraphics[scale=0.11]{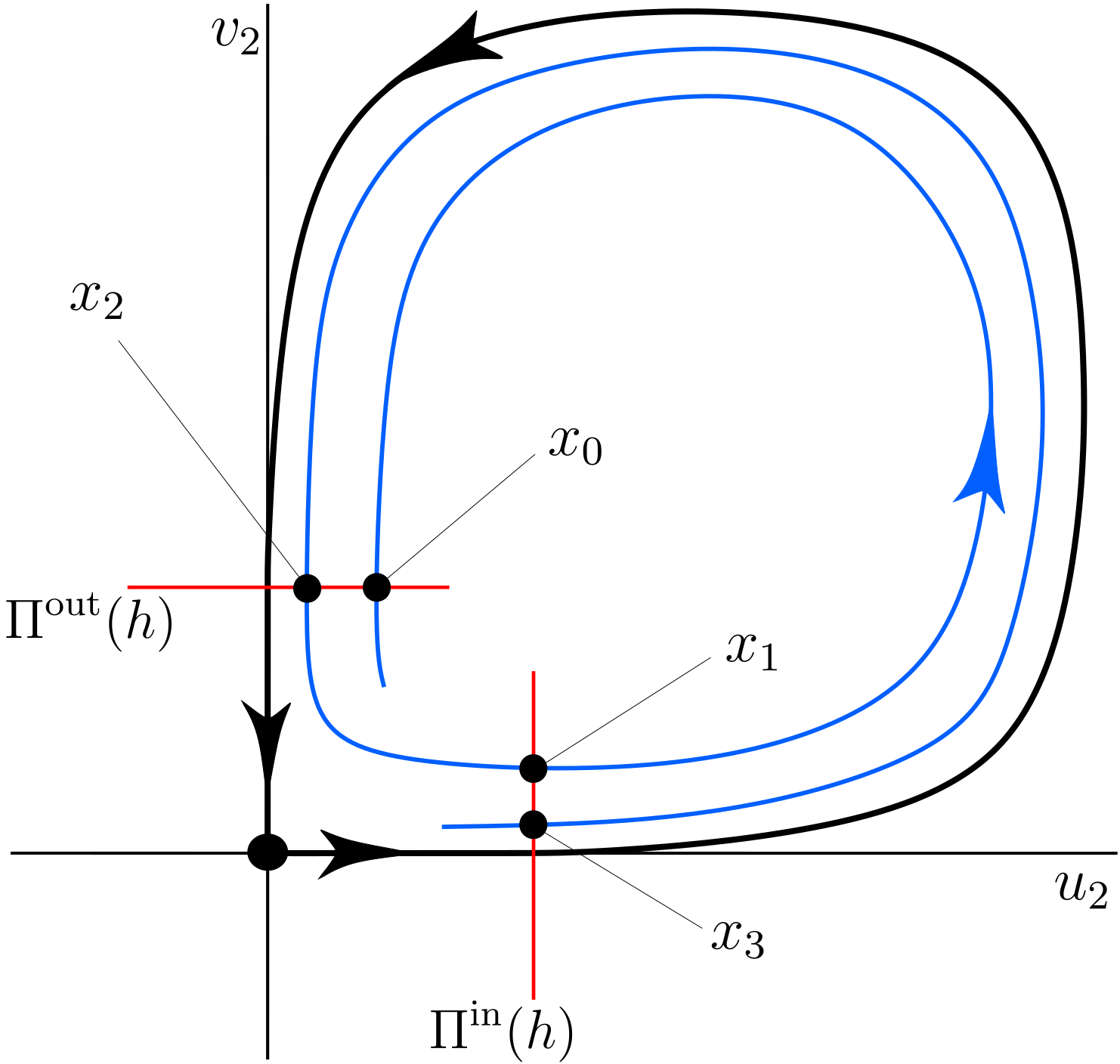}
\caption{System II}
\end{subfigure}
\begin{subfigure}{0.32\textwidth}
\centering
\includegraphics[scale=0.11]{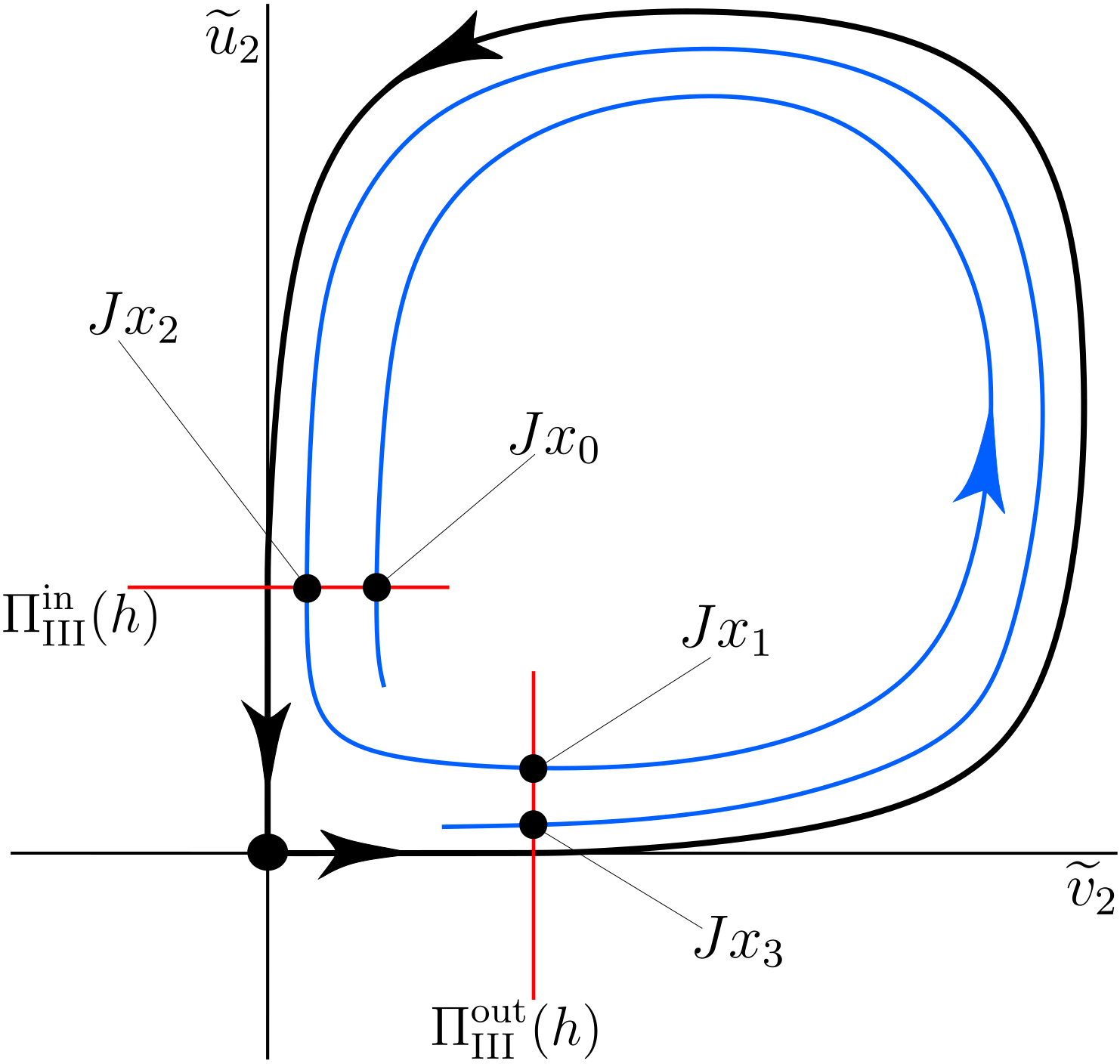}
\caption{System III}
\end{subfigure}
\caption{\small System I is the system in the normal form. We obtain system II from system I by a reversion of time. We obtain system III from system II by changing the stable and unstable components.}
\label{Reversed_time_systems}
\end{figure}

System III is of the form of the normal form that we started with, where all the conditions we assumed for the normal form also hold for system III. This property resembles time-reversibility; however, note that our system is not necessarily time-reversible. The Poincar\'e map of system III maps $Jx_{2}$ to $Jx_{0}$. The global map of system III is \(J \left(T^{\text{glo}}\right)^{-1} J^{-1}\), where \(T^{\text{glo}}\) is the global map of system I. Thus, the differential of this map at $(0,0)$ is
\begin{equation*}
J \left(D T^{glo}\left(M^{s}\right)\right)^{-1} J^{-1} = \frac{1}{a(h) d(h) - b(h) c(h)}\left(\begin{matrix}
a(h) & -c(h)\\
-b(h) & d(h)
\end{matrix}\right).
\end{equation*}
Following the first part of Lemma \ref{Lem11u67tt6} that we have just proved above, if $Jx_{2}\in D_{2}:= \{(\widetilde{u}_{1}, \widetilde{v}_{1})\in\Pi^{\mathrm{in}}_{\mathrm{III}}(h):\,\, m^{-1}\lvert \widetilde{u}_{1}\rvert \leq \lvert \widetilde{v}_{1} \rvert\}$, then $J x_{0}\in D_{2}$, and $Jx_{2} = o\left(\|Jx_{0}\|\right)$ as $(\epsilon, h_{0}) \rightarrow (0,0)$. Moreover, if we write $Jx_{0} = (\widetilde{u}_{1},\widetilde{v}_{1})$, then as $x_{2}\rightarrow (0,0)$, we have $\frac{\widetilde{v}_{1}}{\widetilde{u}_{1}}\rightarrow \frac{-d(h)}{c(h)}$. This implies that for system II, if $x_{2} \in R:= J^{-1}D_{2} = \{(u_{1}, v_{1})\in\Pi^{\mathrm{out}}(h):\,\, \lvert v_{1}\rvert \leq m\lvert u_{1}\rvert\}$, then $x_{0} \in R$, and $x_{2} = o\left(\|x_{0}\|\right)$ as $(\epsilon, h_{0}) \rightarrow (0,0)$. Furthermore, if we write $x_{0} = (u_{1}, v_{1})$, then $\frac{v_{1}}{u_{1}} \rightarrow \frac{-c(h)}{d(h)}$ as $x_{2}\rightarrow (0,0)$.

Let $\widetilde{T}$ be the Poincar\'e map of system II (note that the homoclinic orbit in system II approaches $O$ along $v_{2}$-axis, and so $\widetilde{T}$ is defined on a subset of $\Pi^{\mathrm{out}}(h)$ into $\Pi^{\mathrm{out}}(h)$). In particular, we have $\widetilde{T}(x_{2}) = x_{0}$. Consider the map $\widehat{T}:= T^{\mathrm{glo}}\circ\widetilde{T}\circ\left(T^{\mathrm{glo}}\right)^{-1}$. Note that $\widehat{T}(x_{3}) = x_{1}$. Again, if $x_{3} = T^{\mathrm{glo}}(x_{2})\in T^{\mathrm{glo}}(R)$, then $x_{1} = \widehat{T}(x_{3})\in T^{\mathrm{glo}}(R)$, and $x_{3} = o\left(\|x_{1}\|\right)$ as $(\epsilon, h_{0}) \rightarrow (0,0)$. Furthermore, if we write $x_{1} = (u_{1}, v_{1})$, then $\frac{v_{1}}{u_{1}} \rightarrow 0$ as $x_{3}\rightarrow (0,0)$. The latter is derived from
\begin{equation*}
D T^{\text{glo}}\left(0,0\right)\left(
\begin{matrix}
d(h)\\
-c(h)
\end{matrix}\right) =
\left(\begin{matrix}
a(h) & c(h)\\
b(h) & d(h)
\end{matrix}\right) \left(
\begin{matrix}
d(h)\\
-c(h)
\end{matrix}\right) = \left(
\begin{matrix}
a(h) d(h) - b(h) c(h)\\
0
\end{matrix}\right)
\end{equation*}
and the fact that $a(h) d(h) - b(h) c(h)\neq 0$ which holds since $T^{\mathrm{glo}}$ is a diffeomorphism.

Since the point $x_{3}$ is arbitrary, we have that $\widehat{T}$ and $T^{-1}$, where $T$ is the Poincar\'e map of system I, coincide. Thus, we are done with the proof of the second part of Lemma \ref{Lem11u67tt6} once we show $\mathcal{D}_{1}\cap T(\mathcal{D}_{1}) \subseteq T^{\mathrm{glo}}(R)$. To show this, note that if $k$ is the slope of a point in $\mathcal{D}_{1}$, then $\lvert k\rvert < m^{-1}$. Let $l_{k}$ be the straight line with the slope $k$. The differential of $(T^{\mathrm{glo}})^{-1}$ maps $l_{k}$ to the straight line on $\Pi^{\mathrm{out}}(h)$ that passes through $(0,0)$ and its slope is $\frac{d(h) - k b(h)}{-c(h) + k a(h)}$. If $m$ is large enough, this line lies in $R$. Thus, $\mathcal{D}_{1}\cap T(\mathcal{D}_{1}) \subseteq T^{\mathrm{glo}}(R)$ if $\epsilon > 0$ is chosen sufficiently small. This ends the proof of the second part of Lemma \ref{Lem11u67tt6}.
\end{proof}

We are now ready to prove Theorem \ref{Thm_Single_Loop}. Assume $h < 0$. The periodic orbit $L_h$ intersects $\Pi^{\mathrm{in}}(h)$ at $M^{\mathrm{in}}(h) = (0,0)$, and so, the domain of $T$ contains at least a small neighborhood of $M^{\mathrm{in}}(h)$ in $\Pi^{\mathrm{in}}(h)$ (see Lemma \ref{Lem76tgitiewr}). Let $A := DT(0,0)$ with eigenvalues $\alpha$ and $\beta$. Without loss of generality, suppose $\lvert\alpha\rvert \leq \lvert\beta\rvert$. We claim $\lvert\alpha\rvert < 1 < \lvert\beta\rvert$. To prove this claim, write $T(u_{10}, v_{10}) = A \left(\begin{smallmatrix}
u_{10}\\ v_{10}
\end{smallmatrix}\right) + o\left(\|(u_{10}, v_{10})\|\right)$. By Lemma \ref{Lem11u67tt6}, we have $\frac{\overline{v}_{10}}{\overline{u}_{10}} \rightarrow \frac{d(h)}{b(h)}$ as $(u_{10}, v_{10})\rightarrow (0,0)$. Note that, since $b$ and $d$ are smooth functions of $h$, the assumption $bd\neq 0$ of Theorem \ref{Thm_Single_Loop} implies that $\frac{d(h)}{b(h)}$ is well-defined and nonzero. Thus, if $m > 1$ is large enough, we have that $A$ preserves and expands the cone
\begin{equation*}
C^{u} := \{x:= (u_{10}, v_{10})\in \mathbb{R}^{2}:\quad m^{-1}\lvert u_{10}\rvert \leq \lvert v_{10}\rvert\}.
\end{equation*}
Therefore, $\|Ax\|> K\|x\|$ for some $K > 1$ and all $x\in C^{u}$. Thus, $\|A^{n}\| > K^{n}$ for all $n\in\mathbb{N}$. However, it is known that $\sqrt[n]{\|A^{n}\|}$ converges to the spectral radius of $A$, i.e. $\lvert \beta\rvert$, as $n\rightarrow\infty$, which implies $\lvert \beta\rvert\geq K > 1$. This proves $\lvert \beta\rvert > 1$. Applying the same procedure on $T^{-1}$, the matrix $A^{-1}$ and the cone $C^{s} := \{(u_{10}, v_{10})\in \mathbb{R}^{2}:\,\, \lvert v_{10}\rvert \leq m^{-1}\lvert u_{10}\rvert\}$, we obtain $\lvert \alpha\rvert < 1$. This proves the claim.

Having $\lvert\alpha\rvert < 1 < \lvert\beta\rvert$ implies that the fixed point $(0,0)$ of $T$ possesses a one-dimensional local unstable invariant manifold and a one-dimensional local stable invariant manifold. Extend these manifolds to the whole ball $B_{\epsilon}$ and call them $\Lambda^{u}$ and $\Lambda^{s}$, respectively. Note that $\Lambda^{u}$ and $\Lambda^{s}$ lie in $\mathcal{D}_2$ and $\mathcal{D}_1$, respectively, and are tangent to the straight line whose slope is $\frac{d(h)}{b(h)}$ and the horizontal axis $\{v_{10} = 0\}$, respectively, at $(0,0)$. On the other hand, by Corollary \ref{Cor24g76ftyfduwq}, the points in $\mathcal{D}$ whose backward orbits remain in $\mathcal{D}$ must lie entirely in $\mathcal{D}_{2}$ and their backward orbits converge to $(0,0)$. This means that such points must lie on $\Lambda^u$. Thus, $\Lambda^u$ is exactly the set of the points that never leave $\mathcal{D}$ under the backward iterations of $T$. The proof of the case of forward iterations is analogous. By applying the flow of the system on $\Lambda^{s}$ and $\Lambda^{u}$, we obtain Theorem \ref{Thm_Single_Loop} for the case $h<0$.

Now consider the case $h>0$ and suppose that the backward iterations of a point in $\mathcal{D}$ remain in $\mathcal{D}$. By Corollary \ref{Cor24g76ftyfduwq}, the backward iterations must converge to $(0,0)$. However, by Lemma \ref{Lem76tgitiewr}, the domain $\mathcal{D}$ has no point of a small open neighborhood of $(0,0)$. This contradicts the assumption that the backward iterations of the point remain entirely in $\mathcal{D}$. Therefore, all points in $\mathcal{D}$ leave $\mathcal{D}$ after a finite number of the backward iterations of $T$. The proof of the case of the forward iterations is analogous.

To finalize the proof, we remark on the choice of $\epsilon$. The statement of Lemma \ref{Lem11u67tt6} holds for a sufficiently small $\epsilon$ and all $h\neq 0$ such that $\lvert h\rvert \leq h_{0}$, where $h_{0}$ is sufficiently small. On the other hand, the main theorem of \cite{Bakrani2022JDE} for the dynamics near $\Gamma$ in the level set $\{H=0\}$ holds for some sufficiently small $\epsilon_{0}$. To claim that we describe the dynamics near $\Gamma$ in an open neighborhood $U$ of $\Gamma\cup \{O\}$ in $\mathbb{R}^{4}$, it is sufficient to choose $\epsilon < \epsilon_{0}$. This ends the proof of Theorem \ref{Thm_Single_Loop}.

\section{The Poincar\'e maps and the proof of Lemma \ref{Lem568g67eywfg76}}\label{Technical_Section}

In this section, we explore the Poincar\'e maps for all the four cases $\lambda_{1} = \lambda_{2}$, $\lambda_{1} < \lambda_{2} < 2\lambda_{1}$, $\lambda_{2} = 2\lambda_{1}$ and $2\lambda_{1} < \lambda_{2}$, and prove Lemma \ref{Lem568g67eywfg76}. We start with the case $\lambda_{1} = \lambda_{2}$.

\begin{figure}[h]
\centering
\begin{subfigure}{0.4\textwidth}
\centering
\includegraphics[scale=0.18]{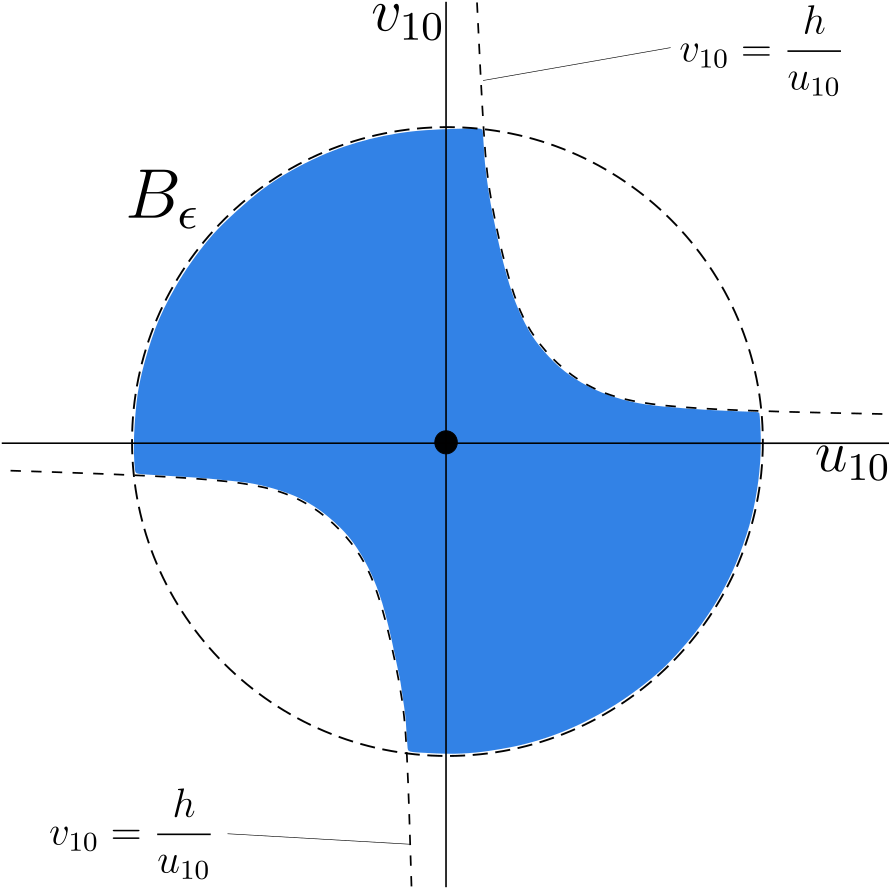}
\caption{\small The case $h < 0$.}
\end{subfigure}
\begin{subfigure}{0.4\textwidth}
\centering
\includegraphics[scale=0.18]{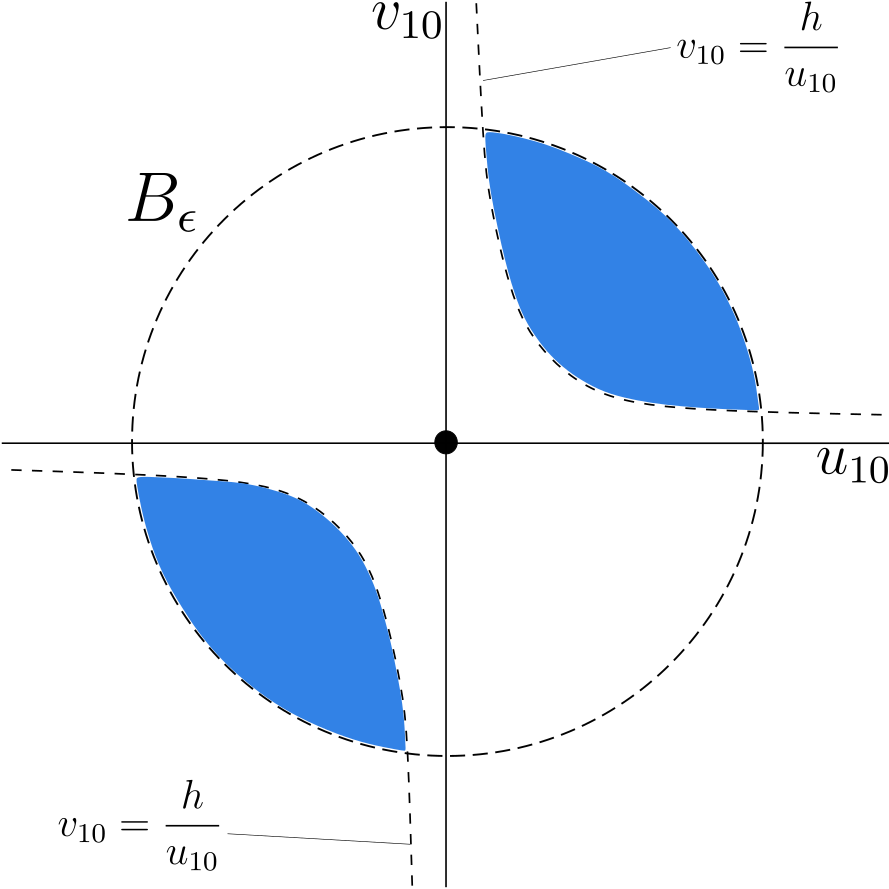}
\caption{\small The case $h > 0$.}
\end{subfigure}
\caption{\small The domain $\mathcal{D}\subset \Pi^{\mathrm{in}}(h)$ of the Poincar\'e map for system (\ref{eq74wt366cuw6gt5uw6v01092}), i.e. the case $\lambda_{1} = \lambda_{2}$, is a subset of the blue region. Indeed, the domain $\mathcal{D}$ is the set of those points in the blue region for which $\|(u_{1\tau}, v_{1\tau})\| \leq \epsilon$. By Proposition \ref{Prop_Domain_is_empty}, if $\epsilon>0$ is sufficiently small, then $\mathcal{D} = \emptyset$.}
\label{Figure_Domain_gamma_equal_1}
\end{figure}

Consider system (\ref{eq74wt366cuw6gt5uw6v01092}). The trajectories near the origin are estimated by (\ref{eq65f3o8fgrcdc}). Evaluating the first two equations of (\ref{eq65f3o8fgrcdc}) at $t = 0$ and $u_{20} = \delta$, and the last two equations at $t=\tau$ and $v_{2\tau} = \delta$ gives
\begin{align}
u_{1\tau} &= e^{-\lambda \tau} u_{10}\left[1 + O\left(\delta\right)\right] + e^{-\lambda\tau} O\left(\lvert v_{1\tau}\rvert\delta\right),\label{eq77m111}\\
u_{2\tau} &= e^{-\lambda \tau} \delta\left[1 + O\left(\delta\right)\right],\nonumber\\
v_{10} &= e^{-\lambda \tau} v_{1\tau}\left[1 + O\left(\delta\right)\right] + e^{-\lambda\tau} O\left(\lvert u_{10}\rvert\delta\right),\label{eq77m333}\\
v_{20} &= e^{-\lambda \tau} \delta\left[1 + O\left(\delta\right)\right].\label{eq77m444}
\end{align}
Relation (\ref{eq77m333}) gives
\begin{equation}\label{eq7rt6t111}
v_{1\tau} = e^{\lambda \tau} v_{10}\left[1 + O\left(\delta\right)\right] + O\left(\lvert u_{10}\rvert\delta\right).
\end{equation}
Substituting this relation into (\ref{eq77m111}) gives
\begin{equation}\label{eq7rt6t222}
u_{1\tau} = e^{-\lambda \tau} u_{10}\left[1 + O\left(\delta\right)\right] + O\left(\lvert v_{10}\rvert\delta\right).
\end{equation}
The point $(u_{10}, \delta, v_{10}, v_{20})\in\Pi^{\mathrm{in}}(h)$ lies in the level set $\{H=h\}$, where $H$ is given by (\ref{eq98n9uo8wbr8bau4b819ub92q1u8}). Thus,
\begin{equation}\label{eq16ruytg6}
v_{20} = \delta^{-1} (u_{10} v_{10} - h).
\end{equation}
Relation (\ref{eq77m444}) gives $v_{20} > 0$, and so by (\ref{eq16ruytg6}), we have $u_{10} v_{10} > h$. This relation implies
\begin{equation*}
\mathcal{D} \subset \{(u_{10}, v_{10}): \quad u_{10} v_{10} > h, \quad  \|(u_{10}, v_{10})\|\leq \epsilon \quad \mathrm{and}\quad \|(u_{1\tau}, v_{1\tau})\|\leq \epsilon\}
\end{equation*}
for an $\epsilon > 0$ (see Figure \ref{Figure_Domain_gamma_equal_1}).

Relations (\ref{eq77m444}) and (\ref{eq16ruytg6}) give $e^{-\lambda\tau} = \delta^{-2}\left[u_{10} v_{10} - h\right]\left[1 + O\left(\delta\right)\right]$, and so
\begin{equation}\label{eq65yy222}
e^{\lambda\tau} = \delta^{2} \left[1 + O\left(\delta\right)\right](u_{10} v_{10} - h)^{-1}.
\end{equation}
Substituting these two relations into (\ref{eq7rt6t111}) and (\ref{eq7rt6t222}) gives an explicit formula for $T^{\mathrm{loc}}$ in terms of $(u_{10}, v_{10})$, and so by (\ref{eq_composition_of_global_and_local_maps}) and (\ref{eq_Taylor_global_map}), we obtain an explicit formula for the Poincar\'e map $T$.

\begin{proof}[Proof of Lemma \ref{Lem568g67eywfg76} for the case $\lambda_{1} = \lambda_{2}$]
Consider $(u_{10}, v_{10}) \in \mathcal{D}_2$. Thus, $u_{10} = O(\lvert v_{10}\rvert)$. Substituting this into (\ref{eq7rt6t111}) gives
\begin{equation}\label{eq872i7ritefuwud}
v_{1\tau} = e^{\lambda\tau} v_{10}\left[1 + O(\delta)\right].
\end{equation}
Then, by (\ref{eq77m111}) and taking $\delta > 0$ small enough such that $\lvert O(\delta)\rvert < 1$, we have
\begin{equation}\label{eq87yoiuhkf6}
\lvert u_{1\tau}\rvert \lvert v_{1\tau}\rvert^{-1} \leq e^{-2\lambda\tau} \lvert u_{10}\rvert \lvert v_{10}\rvert^{-1}  \left[1 + O\left(\delta\right)\right] + e^{-\lambda\tau} O\left(\delta\right) \leq m e^{-2\lambda \tau}  + e^{-\lambda\tau} O\left(\delta\right).
\end{equation}
Therefore, as $(h_{0}, \epsilon)\rightarrow (0,0)$, we have $\tau \rightarrow \infty$, and so $u_{1\tau} = o(\lvert v_{1\tau}\rvert)$, as desired.
\end{proof}

In the proof of Theorem \ref{Thm_Single_Loop}, we showed that when $h>0$, all points in $\mathcal{D}$ leave $\mathcal{D}$ after a finite number of forward and backward iterations. For the case of $\lambda_{1} = \lambda_{2}$, we can actually prove a stronger statement as follows.
\begin{myprop}\label{Prop_Domain_is_empty}
Consider system (\ref{eq74wt366cuw6gt5uw6v01092}). There exists $h_{0} > 0$ for which there exists a sufficiently small $\epsilon > 0$ such that for all $0 < h < h_{0}$, we have $\mathcal{D} = \emptyset$.
\end{myprop}

\begin{proof}
Let $(u_{10}, v_{10})\in \mathcal{D}_{2}$. Then, (\ref{eq65yy222}) and (\ref{eq872i7ritefuwud}) give $v_{1\tau} = \delta^{2} v_{10} \left[1 + O\left(\delta\right)\right] (u_{10} v_{10} - h)^{-1}$. By virtue of this relation, and taking $u_{10} v_{10} - h > 0$ and $h > 0$ into account, we obtain
\begin{equation*}
\epsilon \geq \|(u_{1\tau}, v_{1\tau})\| \geq
\lvert v_{1\tau}\rvert = \frac{\delta^2 \lvert v_{10}\rvert \left[1 + O(\delta)\right]}{u_{10} v_{10} - h} \geq
\frac{\delta^2 \lvert v_{10}\rvert \left[1 + O(\delta)\right]}{u_{10} v_{10}} \geq
\frac{\delta^2 \left[1 + O(\delta)\right]}{\lvert u_{10}\rvert} \geq
\frac{\delta^2 \left[1 + O(\delta)\right]}{\epsilon},
\end{equation*}
and so $\epsilon^{2} \geq \delta^{2} [1 + O(\delta)]$. However, this cannot happen if $\epsilon$ is sufficiently small. Therefore, $\mathcal{D}_2 = \emptyset$. By a reversion of time as in the proof of Lemma \ref{Lem11u67tt6}, we also have $\mathcal{D}_1 = \emptyset$. Thus, $\mathcal{D} = \mathcal{D}_{1}\cup \mathcal{D}_{2} = \emptyset$.
\end{proof}

We now study the Poincar\'e map $T$ for the case $\lambda_{1} < \lambda_{2} < 2\lambda_{1}$. Consider system (\ref{eq19000}). The trajectories near the origin are estimated by (\ref{eqiomljd784987qpqpla874}). Evaluating the first two equations of (\ref{eqiomljd784987qpqpla874}) at $t = 0$ and $u_{20} = \delta$, and the last two equations at $t=\tau$ and $v_{2\tau} = \delta$ gives
\begin{align}
u_{1\tau} =& e^{-\lambda_{1}\tau}u_{10}\left[1 + O\left(\delta\right)\right] + e^{- \lambda_{2}\tau} O\left(\delta \lvert v_{1\tau}\rvert\right),\label{eq17e65er1ytf2ue56f2}\\
u_{2\tau} =& e^{-\lambda_{2}\tau}\delta\left[1 + O\left(\delta\right)\right],\nonumber\\
v_{10} =& e^{-\lambda_{1}\tau} v_{1\tau} \left[1 + O\left(\delta\right)\right] + e^{-\lambda_{2}\tau} O\left(\delta \lvert u_{10}\rvert\right),\label{eq3hukigyuirgf33}\\
v_{20} =& e^{-\lambda_{2}\tau}\delta\left[1 + O\left(\delta\right)\right].\label{eq6739736789927600657}
\end{align}
From (\ref{eq3hukigyuirgf33}), we obtain
\begin{equation}\label{eqyusolkuoplkmh23ui8hd7837}
v_{1\tau} = e^{\lambda_{1}\tau} v_{10} \left[1 + O\left(\delta\right)\right] + e^{\left(\lambda_{1}-\lambda_{2}\right)\tau} O\left(\delta \lvert u_{10}\rvert\right).
\end{equation}
By (\ref{eq6739736789927600657}), we have $v_{20} > 0$. On the other hand, since $(u_{10}, \delta, v_{10}, v_{20})\in\Pi^{\mathrm{in}}(h)$  lies in the level set $\{H=h\}$, where $H$ is given by (\ref{eq22000yuio}), we have that
\begin{equation}\label{eq7yoi7hqi6gfutyvjhby}
v_{20} = \gamma \delta^{-1} u_{10}v_{10}\left[1+o(1)\right] - h \delta^{-1} \left[1+o(1)\right].
\end{equation}
Thus, considering Figure \ref{Figure_Domain_gamma_equal_1} again, the domain here is a subset of the blue region if the curves $v_{10} = \frac{h}{u_{10}}$ in that figure are replaced by $v_{10} = \frac{h}{\gamma u_{10}}\left[1+o(1)\right]$.

Substituting (\ref{eq7yoi7hqi6gfutyvjhby}) into (\ref{eq6739736789927600657}) gives
\begin{equation*}
e^{-\lambda_{2}\tau} = \gamma\delta^{-2} u_{10}v_{10}\left[1 + O\left(\delta\right)\right] - h\delta^{-2} \left[1 + O\left(\delta\right)\right].
\end{equation*}
Having this expression for $e^{-\lambda_{2}\tau}$, we can also express $e^{-\lambda_{1}\tau}$ in terms of $(u_{10}, v_{10})$ using the identity $e^{-\lambda_{1}\tau} = \left(e^{-\lambda_{2}\tau}\right)^{\gamma}$. Substituting these two expressions into (\ref{eq17e65er1ytf2ue56f2}) and (\ref{eqyusolkuoplkmh23ui8hd7837}) gives an explicit formula for the local map $T^{\mathrm{loc}}$, and so the Poincar\'e map $T$.

\begin{proof}[Proof of Lemma \ref{Lem568g67eywfg76} for the case $\lambda_{1} < \lambda_{2} < 2\lambda_{1}$]
Let $(u_{10}, v_{10}) \in \mathcal{D}_{2}$. Therefore, $u_{10} = O(\lvert v_{10}\rvert)$, and so by (\ref{eqyusolkuoplkmh23ui8hd7837}), we have $v_{1\tau} = e^{\lambda_{1}\tau} v_{10} \left[1 + O\left(\delta\right)\right]$. Thus, by (\ref{eq17e65er1ytf2ue56f2}) and taking $\delta > 0$ small enough so that $\lvert O(\delta)\rvert < 1$, we obtain
\begin{equation*}
\lvert u_{1\tau}\rvert \lvert v_{1\tau}\rvert^{-1} \leq e^{-2\lambda_{1}\tau} \lvert u_{10}\rvert \lvert v_{10}\rvert^{-1} \left[1 + O\left(\delta\right)\right] + e^{-\lambda_{2}\tau} O\left(\delta\right)
\leq 2m e^{-2\lambda_{1}\tau} + e^{-\lambda_{2}\tau} \leq 3m e^{-\lambda_{2}\tau}.
\end{equation*}
Therefore, as $(h_{0}, \epsilon)\rightarrow (0,0)$, we have $\tau \rightarrow \infty$, and so $u_{1\tau} = o(\lvert v_{1\tau}\rvert)$, as desired.
\end{proof}

We now analyze the case $\lambda_{2}= 2\lambda_{1}$. Consider system (\ref{eq19000}). The trajectories near the origin are given by (\ref{equ13rhuef8yreiosueygs}). Evaluating the first two equations of (\ref{equ13rhuef8yreiosueygs}) at $t = 0$ and $u_{20} = \delta$, and the last two equations at $t=\tau$ and $v_{2\tau} = \delta$ gives
\begin{align}
u_{1\tau} =& e^{-\lambda \tau} u_{10}\left[1 + O(\delta)\right] + \tau e^{-2\lambda\tau} O\left(\lvert v_{1\tau}\rvert\delta\right),\label{eq87yquigefgew111}\\
u_{2\tau} =& e^{-2\lambda \tau} \delta\left[1 + O\left(\delta\right)\right] + \tau e^{-2\lambda \tau} O\left(u_{10}^{2}\right),\nonumber\\
v_{10} =& e^{-\lambda\tau} v_{1\tau} \left[1 + O\left(\delta\right)\right] + \tau e^{-2\lambda \tau} O\left(\lvert u_{10}\rvert \delta\right),\label{eq232tfhgdvt333}\\
v_{20} =& e^{-2\lambda\tau} \delta\left[1 + O\left(\delta\right)\right] + \tau e^{-2\lambda\tau} O\left(v_{1\tau}^{2}\right).\label{eq2287eh7e444}
\end{align}
Despite the previous two cases, we cannot obtain that much information about the domain $\mathcal{D}$ as well as an explicit formula for $\tau$ in terms of $(u_{10}, v_{10})$ from (\ref{eq2287eh7e444}). However, we can still prove Lemma \ref{Lem568g67eywfg76}.

\begin{proof}[Proof of Lemma \ref{Lem568g67eywfg76} for the case $\lambda_{2} = 2\lambda_{1}$]
Equation (\ref{eq232tfhgdvt333}) yields
\begin{equation*}
v_{1\tau} = e^{\lambda\tau} v_{10} \left[1 + O\left(\delta\right)\right] + \tau e^{-\lambda \tau} O\left(\lvert u_{10}\rvert \delta\right).
\end{equation*}
In particular, if $(u_{10}, v_{10})\in \mathcal{D}_{2}$, we have $u_{10} = O(\lvert v_{10}\rvert)$, and so by taking into account that $\tau e^{-\tau}\rightarrow 0$ as $\tau\rightarrow\infty$ (this holds as $(h_{0}, \epsilon)\rightarrow (0,0)$), we get $v_{1\tau} = e^{\lambda\tau} v_{10} \left[1 + O\left(\delta\right)\right]$. Thus, by (\ref{eq87yquigefgew111}), we obtain
\begin{equation*}
\lvert u_{1\tau}\rvert \lvert v_{1\tau}\rvert^{-1} \leq e^{-2\lambda\tau} \lvert u_{10}\rvert \lvert v_{10}\rvert^{-1}  \left[1 + O\left(\delta\right)\right] + \tau e^{-2\lambda\tau} O\left(\delta\right)
\leq 2m e^{-2\lambda\tau} + e^{-\lambda\tau} \leq 3m e^{-\lambda\tau}.
\end{equation*}
Therefore, as $(h_{0}, \epsilon)\rightarrow (0,0)$, we have $\tau \rightarrow \infty$, and so $u_{1\tau} = o(\lvert v_{1\tau}\rvert)$, as desired.
\end{proof}

Finally, we study the Poincar\'e map for the case $\lambda_{2} > 2\lambda_{1}$. Consider system (\ref{eq23000}). The trajectories near the origin are estimated by (\ref{equos83km3odnin8b83xbv}). Evaluating the first two equations of (\ref{equos83km3odnin8b83xbv}) at $t = 0$ and $u_{20} = \delta$, and the last two equations at $t=\tau$ and $v_{2\tau} = \delta$ gives
\begin{align}
u_{1\tau} &= e^{-\lambda_{1}\tau}u_{10}\left[1 + O\left(\delta\right)\right] + e^{-2\lambda_{1}\tau} O\left(\delta \lvert v_{1\tau}\rvert\right),\label{eq1r115u5ruytgsw23}\\
u_{2\tau} &= e^{-\lambda_{2}\tau} \delta \left[1 + O\left(\delta\right)\right],\nonumber\\
v_{10} &= e^{-\lambda_{1}\tau}v_{1\tau}\left[1 + O\left(\delta\right)\right] + e^{-2\lambda_{1}\tau} O\left(\delta \lvert u_{10}\rvert\right),\label{eq7198ytyqi6ugtef}\\
v_{20} &= e^{-\lambda_{2} \tau} \delta \left[1 + O\left(\delta\right)\right].\label{eq149tyyh3ehlngt}
\end{align}
In particular, (\ref{eq7198ytyqi6ugtef}) gives
\begin{equation}\label{eq5tutfytyvttcdssw}
v_{1\tau} = e^{\lambda_{1}\tau}v_{10}\left[1 + O\left(\delta\right)\right] + e^{-\lambda_{1}\tau} O\left(\delta \lvert u_{10}\rvert\right).
\end{equation}
The point $(u_{10}, \delta, v_{10}, v_{20})\in\Pi^{\mathrm{in}}(h)$ lies in the level set $\{H=h\}$, where $H$ is given by (\ref{eq12hniuwb7gydcw}). Thus,
\begin{equation}\label{eq2308u7yi7go8yed}
h = \gamma u_{10}v_{10} \left[1 + o\left(1\right)\right] - v_{20}\delta \left[1 + o\left(1\right)\right] + v_{10}^{2} O\left(\delta\right) + v_{20} u_{10}^{2} O\left(1\right).
\end{equation}
Note that $u_{10} = o(1)$ as $\epsilon\rightarrow 0$, and so $v_{20} u_{10}^{2} O\left(1\right)$ is of the order of $v_{20}\delta o(1)$. Thus, (\ref{eq2308u7yi7go8yed}) yields
\begin{equation}\label{eq32jy8rd58tiuy76g}
v_{20} = \delta^{-1} \left(\gamma u_{10}v_{10} \left[1 + o\left(1\right)\right] + v_{10}^{2} O\left(\delta\right) - h \left[1 + o\left(1\right)\right]\right).
\end{equation}
Following (\ref{eq149tyyh3ehlngt}) and (\ref{eq32jy8rd58tiuy76g}), we have
\begin{equation*}
e^{-\lambda_{2}\tau} = \delta^{-2}\left(\gamma u_{10}v_{10} \left[1 + o\left(1\right)\right] + v_{10}^{2} O\left(\delta\right) - h\left[1 + o\left(1\right)\right]\right).
\end{equation*}
This relation gives an explicit formula for $e^{-\lambda_{2}\tau}$ in terms of $(u_{10}, v_{10})$ even though none of the terms $\gamma u_{10}v_{10} \left[1 + o\left(1\right)\right]$ and $v_{10}^{2} O\left(\delta\right)$ is dominant. Then, by the identity $e^{-\lambda_{1}\tau} = \left(e^{-\lambda_{2}\tau}\right)^{\gamma}$, an expression for $e^{-\lambda\tau}$ in terms of $(u_{10}, v_{10})$ is obtained as well, and so by substituting these into (\ref{eq1r115u5ruytgsw23}) and (\ref{eq5tutfytyvttcdssw}), we find an explicit formula for $T^{\mathrm{loc}}$, and so for the Poincar\'e map $T$, too.

It also follows from (\ref{eq149tyyh3ehlngt}) that for the domain $\mathcal{D}$, we must have $v_{20} > 0$. However, we cannot express this relation terms of $(u_{10}, v_{10})$ using (\ref{eq32jy8rd58tiuy76g}) as none of the terms $\gamma u_{10}v_{10} \left[1 + o\left(1\right)\right]$ and $v_{10}^{2} O\left(\delta\right)$ in (\ref{eq32jy8rd58tiuy76g}) is dominant.

\begin{proof}[Proof of Lemma \ref{Lem568g67eywfg76} for the case $\lambda_{2} > 2\lambda_{1}$]
Consider $(u_{10}, v_{10}) \in \mathcal{D}_{2}$. Thus, $u_{10} = O(\lvert v_{10}\rvert)$, and by (\ref{eq5tutfytyvttcdssw}), we have $v_{1\tau} = e^{\lambda_{1}\tau} v_{10} \left[1 + O\left(\delta\right)\right]$. Then, by (\ref{eq1r115u5ruytgsw23}) and taking $\delta > 0$ sufficiently small such that $\lvert O(\delta)\rvert < 1$, we obtain
\begin{equation*}
\lvert u_{1\tau}\rvert \lvert v_{1\tau}\rvert^{-1} \leq e^{-2\lambda_{1}\tau} \lvert u_{10}\rvert \lvert v_{10}\rvert^{-1} \left[1 + O\left(\delta\right)\right] + e^{-2\lambda_{1}\tau} O\left(\delta\right) \leq \left(2m+1\right) e^{-2\lambda_{1}\tau}.
\end{equation*}
Therefore, as $(h_{0}, \epsilon)\rightarrow (0,0)$, we have $\tau \rightarrow \infty$, and so $u_{1\tau} = o(\lvert v_{1\tau}\rvert)$, as desired.
\end{proof}

\section{Proof of Theorem \ref{Thm_Double_Loops}} \label{Sec_Proof_Thm_Double_Loops}

We prove Theorem \ref{Thm_Double_Loops} in this section. Recall the homoclinic orbits $\Gamma_{+}$ and $\Gamma_{-}$. Consider two 3-dimensional cross-sections $\Sigma^{\mathrm{in}}_{\sigma} = \{u_2 = \sigma\delta\}$ and $\Sigma^{\mathrm{out}}_{\sigma} = \{v_2 = \sigma\delta\}$ to the loop $\Gamma_{\sigma}$ for $\sigma \in \{+,-\}$, where $\delta > 0$ is sufficiently small (see Figure \ref{Fig_three_dim_cross_sections_double_homoclinics}). Take a sufficiently small $h_{0} > 0$, and let $\Pi^{\mathrm{in}}_{\sigma}(h)$ and $\Pi^{\mathrm{out}}_{\sigma}(h)$ be the restriction of $\Sigma^{\mathrm{in}}_{\sigma}$ and $\Sigma^{\mathrm{out}}_{\sigma}$ to the level set $\{H = h\}$, respectively, where $\lvert h\rvert < h_{0}$ i.e. \(\Pi^{\mathrm{in}}_{\sigma}(h) = \lbrace u_{2} = \sigma\delta\rbrace \cap \lbrace H = h\rbrace\) and \(\Pi^{\mathrm{out}}(h) = \lbrace v_{2} = \sigma\delta\rbrace \cap \lbrace H = h\rbrace\). With the exact same proof, one can show that Lemma \ref{Coordinates_Section_Lemma} holds for $\Pi^{\mathrm{in}}_{\sigma}(h)$ and $\Pi^{\mathrm{out}}_{\sigma}(h)$ too, i.e. we can choose $(u_{1}, v_{1})$-coordinates on $\Pi^{\mathrm{in}}_{\sigma}(h)$ and $\Pi^{\mathrm{out}}_{\sigma}(h)$ meaning that $(u_{2}, v_{2})$ is uniquely determined by $(u_{1}, v_{1})$ on these sections. We also denote the points $(0,0)$ on $\Pi^{\mathrm{in}}_{\sigma}(h)$ and $\Pi^{\mathrm{out}}_{\sigma}(h)$ by $M^{\mathrm{in}}_{\sigma}(h)$ and $M^{\mathrm{out}}_{\sigma}(h)$, respectively.

Fix an $h\in (-h_{0}, h_{0})$, and a sufficiently small $\epsilon >0$. For $\sigma_{1}, \sigma_{2}\in\{+,-\}$, we denote by $\mathcal{D}^{\sigma_{1} \sigma_{2}}$ the set of all points $(u_{10}, v_{10})\in \Pi_{\sigma_{1}}^{\mathrm{in}}(h)$ such that $\|(u_{10}, v_{10})\| < \epsilon$ and their forward orbits go along $\Gamma_{\sigma_{2}}$ and intersect $\Pi^{\mathrm{out}}_{\sigma_{2}}(h)$ at $(u_{1\tau}, v_{1\tau})$, where $\|(u_{1\tau}, v_{1\tau})\| < \epsilon$. We define the local map $T^{\mathrm{loc}}_{\sigma_{1} \sigma_{2}}$ on $\mathcal{D}^{\sigma_{1} \sigma_{2}}$ by $(u_{10}, v_{10})\mapsto (u_{1\tau}, v_{1\tau})$. We also define $T_{\sigma}^{\mathrm{glo}}$ for $\sigma = \pm$ on the open $\epsilon$-ball centered at $(0,0)$ on $\Pi^{\mathrm{out}}_{\sigma}$. Then, the Poincar\'e map along $\Gamma_{\sigma}$ is given by $T_{\sigma} := T^{\mathrm{glo}}_{\sigma}\circ T^{\mathrm{loc}}_{\sigma\sigma}$.

A point on $\Pi^{\mathrm{in}}_{+}(h)$ or $\Pi^{\mathrm{in}}_{-}(h)$ whose forward or backward orbit lies entirely in $V_{h}$ must intersect $\Pi^{\mathrm{in}}_{+}(h)$ or $\Pi^{\mathrm{in}}_{-}(h)$ infinitely many times. If the forward or backward orbit of this orbit intersects only $\Pi^{\mathrm{in}}_{+}(h)$ (resp. $\Pi^{\mathrm{in}}_{-}(h)$), this means that it lies in a small neighborhood of $\Gamma_{+}$ (resp. $\Gamma_{-}$). These orbits are described by Theorem \ref{Thm_Single_Loop}. The remaining scenario is that such orbits intersect each of $\Pi^{\mathrm{in}}_{+}(h)$ and $\Pi^{\mathrm{in}}_{-}(h)$ infinitely many times. To study such orbits, let $\mathcal{D} := \mathcal{D}^{++} \cup \mathcal{D}^{+-} \cup \mathcal{D}^{-+} \cup \mathcal{D}^{--}$. Analogous to (\ref{eqyujnj68999jhhhgg34}), we define the sets $Y_1$ and $Y_2$ for the case of two homoclinic orbits as
\begin{equation*}
\begin{gathered}
Y_{1} := \left\{\left(u_{10}, v_{10}\right)\in \Pi^{\mathrm{in}}_{+}(h) \cup \Pi^{\mathrm{in}}_{-}(h):\,\, \|(u_{10}, v_{10})\| < \epsilon\quad \mathrm{and} \quad \lvert v_{10}\rvert < m^{-1}\lvert u_{10}\rvert\right\} \,\,\mathrm{and}\\
Y_{2} := \left\{\left(u_{10}, v_{10}\right)\in \Pi^{\mathrm{in}}_{+}(h) \cup \Pi^{\mathrm{in}}_{-}(h):\,\, \|(u_{10}, v_{10})\| < \epsilon \quad \mathrm{and} \quad m^{-1}\lvert u_{10}\rvert \leq \lvert v_{10}\rvert\right\}.
\end{gathered}
\end{equation*}
We further define $\mathcal{D}_{i} := \mathcal{D} \cap Y_{i}$. The following lemma is an analog of Lemma \ref{Lem11u67tt6}. It describes the forward dynamics of the points in $\mathcal{D}_{2}$ and the backward dynamics of the points in $\mathcal{D}_1$.

\begin{mylem}\label{Lem563dtyrd5tr}
Reduce system (\ref{eq300}) near the origin to one of the normal forms (\ref{eq74wt366cuw6gt5uw6v01092}), (\ref{eq19000}) or (\ref{eq23000}). Take $(u_{10}, v_{10})\in \mathcal{D}$ and let $(\overline{u}_{10}, \overline{v}_{10})$ be the first intersection point of the forward orbit of $(u_{10}, v_{10})$ and $\Pi^{\mathrm{in}}_{+}(h) \cup \Pi^{\mathrm{in}}_{-}(h)$. The following holds.
\begin{enumerate}[(i)]
\item If $(u_{10}, v_{10})\in \mathcal{D}_2$ such that $\|(\overline{u}_{10}, \overline{v}_{10})\| < \epsilon$, then $(\overline{u}_{10}, \overline{v}_{10}) \in Y_{2}$. Furthermore, $(u_{10}, v_{10}) = o\left(\|(\overline{u}_{10}, \overline{v}_{10})\|\right)$ as $(\epsilon, h_{0})\rightarrow (0,0)$. In addition, for $\sigma = \pm$, as $(u_{10}, v_{10})\rightarrow (0,0)\in \Pi^{\mathrm{in}}_{\sigma}(h)$, we have that $(\overline{u}_{10}, \overline{v}_{10}) \rightarrow (0,0) \in \Pi^{\mathrm{in}}_{\sigma}(h)$ and $\frac{\overline{v}_{10}}{\overline{u}_{10}} \rightarrow \frac{d_{\sigma}(h)}{b_{\sigma}(h)}$.

\item If $(\overline{u}_{10}, \overline{v}_{10}) \in Y_{1}$ for some $(u_{10}, v_{10}) \in \mathcal{D}$, then $(u_{10}, v_{10}) \in \mathcal{D}_{1}$. Furthermore, $(\overline{u}_{10}, \overline{v}_{10}) = o\left(\|(u_{10}, v_{10})\|\right)$ as $(\epsilon, h_{0})\rightarrow (0,0)$. In addition, for $\sigma = \pm$, as $(\overline{u}_{10}, \overline{v}_{10})\rightarrow (0,0) \in \Pi^{\mathrm{in}}_{\sigma}(h)$, we have that $(u_{10}, v_{10})\rightarrow (0,0) \in \Pi^{\mathrm{in}}_{\sigma}(h)$ and $\frac{v_{10}}{u_{10}} \rightarrow 0$.
\end{enumerate}
\end{mylem}

Once we prove a technical statement analogous to Lemma \ref{Lem568g67eywfg76} for the case of double homoclinics, the proof of Lemma \ref{Lem563dtyrd5tr} becomes similar to the proof of Lemma \ref{Lem11u67tt6}. One thing to be taken into account is that for the case of double homoclinics, when $(u_{10}, v_{10}) \in \Pi^{\mathrm{in}}_{\sigma_{1}}(h)$, $(u_{1\tau}, v_{1\tau}) \in \Pi^{\mathrm{out}}_{\sigma_{2}}(h)$ and $(\overline{u}_{10}, \overline{v}_{10}) \in \Pi^{\mathrm{in}}_{\sigma_{2}}(h)$, then despite the case of a single homoclinic, we do not necessarily have $\sigma_{1} = \sigma_{2}$. In this case, if $(u_{10}, v_{10})\in \mathcal{D}_2$ and $\|(\overline{u}_{10}, \overline{v}_{10})\| < \epsilon$, then $(\overline{u}_{10}, \overline{v}_{10})$ is close to either the straight line with the slope $\frac{d_{+}(h)}{b_{+}(h)}$ (the case $(\overline{u}_{10}, \overline{v}_{10}) \in \Pi^{\mathrm{in}}_{+}(h)$) or the straight line with the slope $\frac{d_{-}(h)}{b_{-}(h)}$ (the case $(\overline{u}_{10}, \overline{v}_{10}) \in \Pi^{\mathrm{in}}_{-}(h)$). In any case, once $m>1$ is chosen large enough, we have $(\overline{u}_{10}, \overline{v}_{10})\in Y_2$. To finalize the proof of Lemma \ref{Lem563dtyrd5tr}, we state the following statement that is an analog of Lemma \ref{Lem568g67eywfg76}.

\begin{mylem}
Consider the system in its normal form. Take $(u_{10}, v_{10}) \in \mathcal{D}_{2}$, and let $(u_{1\tau}, v_{1\tau})$ be the first intersection point of the forward orbit of $(u_{10}, v_{10})$ and $\Pi^{\mathrm{out}}_{+}(h) \cup \Pi^{\mathrm{out}}_{-}(h)$. Then, as $(h_{0}, \epsilon) \rightarrow (0,0)$, we have $v_{1\tau} = e^{\lambda_{1}\tau} v_{10}\left[1 + O\left(\delta\right)\right]$ and $u_{1\tau} = o\left(\lvert v_{1\tau}\rvert\right)$.
\end{mylem}

\begin{proof}
Suppose $(u_{10}, v_{10}) \in \Pi^{\mathrm{in}}_{\sigma_{1}}(h)$ and $(u_{1\tau}, v_{1\tau}) \in \Pi^{\mathrm{out}}_{\sigma_{2}}(h)$. The case $\sigma_{1} = \sigma_{2}$ directly follows from Lemma \ref{Lem568g67eywfg76}. The proof of the case $\sigma_{1} \neq \sigma_{2}$ is also similar to the proof of Lemma \ref{Lem568g67eywfg76}. Indeed, one needs to check the scenarios $\lambda_{1} = \lambda_{2}$, $\lambda_{1} < \lambda_{2} < 2\lambda_{1}$, $\lambda_{2} = 2\lambda_{1}$ and $2\lambda_{1} < \lambda_{2}$ one by one. Here, we prove the lemma for the case $\lambda_{1} = \lambda_{2}$; the other three cases are similar.

Assume $\lambda_{1} = \lambda_{2}$. Consider system (\ref{eq74wt366cuw6gt5uw6v01092}). The trajectories near the origin are estimated by (\ref{eq65f3o8fgrcdc}). Evaluating the first two equations of (\ref{eq65f3o8fgrcdc}) at $t = 0$ and $u_{20} = \sigma_{1} \delta$, and the last two equations at $t=\tau$ and $v_{2\tau} = \sigma_{2}\delta$ gives
\begin{align}
u_{1\tau} &= e^{-\lambda \tau} u_{10}\left[1 + O\left(\delta\right)\right] + e^{-\lambda\tau} O\left(\lvert v_{1\tau}\rvert\delta\right),\label{eq1276ty2igd}\\
u_{2\tau} &= \sigma_{1} e^{-\lambda \tau} \delta\left[1 + O\left(\delta\right)\right],\nonumber\\
v_{10} &= e^{-\lambda \tau} v_{1\tau}\left[1 + O\left(\delta\right)\right] + e^{-\lambda\tau} O\left(\lvert u_{10}\rvert\delta\right),\label{eq61tyjgbhoe87}\\
v_{20} &= \sigma_{2} e^{-\lambda \tau} \delta\left[1 + O\left(\delta\right)\right].\nonumber
\end{align}
Relation (\ref{eq61tyjgbhoe87}) gives (\ref{eq7rt6t111}). On the other hand, since $(u_{10}, v_{10})\in \mathcal{D}_{2}$, we have $u_{10} = O(\lvert v_{10}\rvert)$. Substituting $u_{10} = O(\lvert v_{10}\rvert)$ into (\ref{eq7rt6t111}) gives (\ref{eq872i7ritefuwud}), as desired. Then, similar to the proof of Lemma \ref{Lem568g67eywfg76}, relation (\ref{eq1276ty2igd}) implies (\ref{eq87yoiuhkf6}). Thus, as $(h_{0}, \epsilon)\rightarrow (0,0)$, we have $\tau \rightarrow \infty$, and so $u_{1\tau} = o(\lvert v_{1\tau}\rvert)$, as desired.
\end{proof}

We now prove Theorem \ref{Thm_Double_Loops}. Suppose $h < 0$. It follows from Theorem \ref{Thm_Single_Loop} that if the system is restricted to the level set $\{H = h\}$, each of the periodic orbits $L_{h}^{+}$ and $L_{h}^{-}$ are saddle; they possess stable manifolds $W^{s}(L_{h}^{+})$ and $W^{s}(L_{h}^{-})$, and unstable manifolds $W^{u}(L_{h}^{+})$ and $W^{u}(L_{h}^{-})$. Suppose that there exists a point in $V_h$ whose forward orbit remains in $V_h$. By Lemma \ref{Lem563dtyrd5tr}, the intersection points of this forward orbit and $\Pi^{\mathrm{in}}_{+}(h) \cup \Pi^{\mathrm{in}}_{-}(h)$ must belong to $\mathcal{D}_{1} \cup \{(0,0)\}$, and as time passes, these intersection points converge to $(0,0)$. Here, the point $(0,0)$ may correspond to $M^{\mathrm{in}}_{+}(h)$ or $M^{\mathrm{in}}_{-}(h)$. However, due to the smoothness of the system, once the intersection points get close enough to one of the points $M^{\mathrm{in}}_{+}(h)$ or $M^{\mathrm{in}}_{-}(h)$, they remain close to that point forever; in other words, the forward orbit must converge to one of the periodic orbits $L^{+}_{h}$ or $L^{-}_{h}$. Thus, if the forward orbit of a point in $V_h$ remains in $V_h$, that point must belong to $W^{s}(L_{h}^{+}) \cup W^{s}(L_{h}^{-})$. Analogously, one can show that if the backward orbit of a point in $V_h$ remains in $V_h$, that point must belong to $W^{u}(L_{h}^{+}) \cup W^{u}(L_{h}^{-})$. This proves Theorem \ref{Thm_Double_Loops} for the case of negative $h$.

Assume $h>0$. The Poincar\'e map $T$ along $L_h$ defined on a subset of $\mathcal{D}^{+-} \subset \Pi^{\mathrm{in}}_{+}(h)$ is given by $T = T^{\mathrm{glo}}_{+} \circ T^{\mathrm{loc}}_{-+} \circ T^{\mathrm{glo}}_{-} \circ T^{\mathrm{loc}}_{+-}$. Let $A := DT(0,0)$. By Lemma \ref{Lem563dtyrd5tr}, $A$ (resp. $A^{-1}$) preserves and expands the cone $C^{u}$ (resp. $C^{s}$), where $C^{s}$ and $C^{u}$ are the cones introduced in the proof of Theorem \ref{Thm_Single_Loop}. The rest of the proof is the same as the proof of Theorem \ref{Thm_Single_Loop}. It should be noticed that since we defined the Poincar\'e map $T$ on $\Pi^{\mathrm{in}}_{+}(h)$, the curve $\Lambda^{u}$ is tangent to the straight line with the slope $\frac{d_{+}(h)}{b_{+}(h)}$ at $M^{\mathrm{in}}_{+}(h)$, and the curve $T^{\mathrm{glo}}_{-} \circ T^{\mathrm{loc}}_{+-}(\Lambda^{u})$, i.e. the intersection of the unstable manifold of $L_h$ and $\Pi^{\mathrm{in}}_{-}(h)$, is tangent to the straight line with the slope $\frac{d_{-}(h)}{b_{-}(h)}$ at $M^{\mathrm{in}}_{-}(h)$.

\section*{Acknowledgments}

This research was supported by FAPESP (grant no. 2023/04294-0).

%\appendix

%\addtocontents{toc}{\protect\setcounter{tocdepth}{1}}
% This command control what appears on the table of content.

\bibliography{Myreferences}{}
\bibliographystyle{alpha}

\end{document}